\documentclass[12pt]{amsart}
\usepackage{amssymb, amscd,pb-diagram}
\usepackage{verbatim}
\usepackage{hhline}

\usepackage{doublespace}  

\newtheorem{thm}{Theorem}[section]
\newtheorem{prop}[thm]{Proposition}
\newtheorem{lma}[thm]{Lemma}
\newtheorem{cor}[thm]{Corollary}

\newenvironment{defn}{\trivlist
\item[{\bf\hspace{1.88pt}\addtocounter{thm}{1}
Definition~\arabic{section}.\arabic{thm}.}]}{\endtrivlist}

\newcommand{\cal}{\mathcal}
\def\bbC{\mathbb C}
\def\bbN{\mathbb N}

\def\bbZ{\mathbb Z}

\def\sA{{\cal A}}

\def\sC{{\cal C}}

\def\sE{{\cal E}}
\def\sF{\cal F}
\def\sG{{\cal G}}

\def\sK{{\cal K}}

\def\sP{{\cal P}}
\def\sR{{\cal R}}

\def\Ad{{\mathrm{Ad}\,}}

\def\id{{\mathrm{id}}}

\def\Alglim{{\mathrm{Alglim}}}

\def\rad{{\mathrm{rad}\thinspace}}
\def\rank{{\mathrm{rank}\thinspace}}

\def\Lim{{\mathrm{Lim}\thinspace}}
\def\AF{{\mathrm{AF}}}
\def\End{{\mathrm{End}\thinspace}}
\def\Pend{{\mathrm{Pend}\thinspace}}

\def\inf{{\mathrm{inf}\thinspace}}
\def\Sys{{\mathrm{Sys}\thinspace}}

\def\indlimit{{\displaystyle \lim_{\longrightarrow} }}
\def\algindlimit{{\displaystyle \mathop{\mathrm{alg\,lim}}\limits_{\longrightarrow} }}

\def\resultspace{{\vskip 0.075in}}

\def\chisub#1 {{\chi_{\lower2.5pt\hbox{$\scriptstyle #1$}}}}

\def\symbdown#1{\Big\downarrow\rlap{$\vcenter{\hbox{$\scriptstyle #1$}}$}}
\def\symbup#1{\Big\uparrow\rlap{$\vcenter{\hbox{$\scriptstyle #1$}}$}}

\def\symbne#1{\!\!\!\!\nearrow\rlap{$\vcenter{\hbox{$\scriptstyle #1$}}$}}
\def\symbse#1{\!\!\!\!\searrow\rlap{$\vcenter{\hbox{$\scriptstyle #1$}}$}}

\def\symbse#1{\!\!\!\!\searrow\rlap{$\vcenter{\hbox{$\scriptstyle #1$}}$}}

\def\symbup#1{\Big\uparrow\rlap{$\vcenter{\hbox{$\scriptstyle #1$}}$}}
\def\symbdown#1{\Big\downarrow\rlap{$\vcenter{\hbox{$\scriptstyle #1$}}$}}

\def\angb#1{\langle #1 \rangle }
 
\numberwithin{equation}{section}
\begin{document}
\bibliographystyle{amsalpha209}

\begin{spacing}{1.38}


\newcommand{\matrrccc}[6]{\mbox{$ \left[ \begin{array}{ccc}
	    #1&#2&#3\\ #4&#5&#6 \end{array} \right] $}}
\newcommand{\matrccc}[3]{\mbox{$ \left[ \begin{array}{ccc}
	    #1&#2&#3 \end{array} \right] $}}
\newcommand{\matrcccc}[4]{\mbox{$ \left[ \begin{array}{cccc}
	    #1&#2&#3&#4 \end{array} \right] $}}
\newcommand{\matrrrc}[3]{\mbox{$ \left[ \begin{array}{c}
	    #1 \\ #2 \\ #3 \end{array} \right] $}}
\newcommand{\matrrrrc}[4]{\mbox{$ \left[ \begin{array}{c}
	    #1 \\ #2 \\ #3 \\ #4 \end{array} \right] $}}
\newcommand{\matrrrcc}[6]{\mbox{$ \left[ \begin{array}{cc}
	    #1&#2 \\ #3&#4 \\#5&#6 \end{array} \right] $}}
\newcommand{\matrrrrcc}[8]{\mbox{$ \left[ \begin{array}{cc}
	    #1&#2\\#3&#4\\#5&#6\\#7&#8 \end{array} \right] $}}
\newcommand{\matrrrccc}[9]{\mbox{$ \left[ \begin{array}{ccc}                    
            #1&#2&#3\\#4&#5&#6 \\#7&#8&#9 \end{array} \right] $}}
\newcommand{\rrrccc}[9]{\mbox{$  \begin{array}{ccc}                $
            #1&#2&#3\\#4&#5&#6 \\#7&#8&#9 \end{array}  $}}

\newcommand{\mattri}[9]{\mbox{$ \left[ \begin{array}{cccc}
   #1\\#2&#3\\#4&#5&#6 \\#7&#8&#9&1 \end{array} \right] $}}
\newcommand{\mtab}[7]{
\begin{array}{rrr rr}\hline\\[-9pt]
  #1	 & &#3 & 0 & 1 \\[1pt]\hline\\[-9pt]
      #2 &0&  &#4 &#5 \\
	 &1&  &#6 &#7 \\[1pt]\hline
\end{array}  }

\def\angb#1{\langle #1 \rangle }

\newcommand{\matrcc}[2]{\mbox{$ \left[ \begin{array}{cc} #1 & #2 \end{array}
           \right] $}}
\newcommand{\matrrc}[2]{\mbox{$ \left( \begin{array}{c}  #1 \\ #2 \end{array}
           \right) $}}
\newcommand{\matrrcc}[4]{\mbox{$ \left[ \begin{array}{cc} #1&#2 \\ #3&#4
           \end{array} \right] $}}
\newcommand{\IZ}{\mathbb{Z}}

\bibliographystyle{amsalpha209}

\title[Approximately Finitely Acting Operator Algebras]
{Approximately Finitely Acting Operator Algebras}
\date{March, 2000}
\author{Stephen~C. Power}
\address{Dept. of Mathematics \& Statistics \\ Lancaster
 University \\ Lancaster, U.K. LA1 4YF}
\email[]{s.power\@lancaster.ac.uk}

\begin{abstract}
Let $E$ be an  operator algebra on a Hilbert space with
finite-dimensional C*-algebra $C^*(E)$.
A classification is given of the locally finite algebras $A_0 = \algindlimit (A_k,
\phi_k)$ and the operator algebras $A =  \indlimit (A_k,
\phi_k)$ obtained as limits of
direct sums of matrix algebras over $E$ with respect to star-extendible
homomorphisms.
The invariants in the algebraic case consist 
of an additive semigroup, with
scale, which is a right module for the semiring
$~ V_E = Hom_u(E \otimes \sK, E \otimes \sK)
~$
of unitary equivalence classes of star-extendible homomorphisms.
This semigroup is referred to as the dimension module invariant.
In the operator algebra case the invariants consist of a metrized
additive semigroup with scale and a contractive right module $V_E$-action.
Subcategories of algebras determined by restricted classes of
embeddings, such as 1-decomposable embeddings between digraph algebras,
are also classified in terms of simplified dimension modules.
\end{abstract}
\maketitle
 
\tableofcontents

\section{Introduction} 

Approximately finite (AF) C*-algebras are classified in
terms of the scaled $K_0$-group (Elliott \cite{ell-1}). This
perspective subsumed earlier special cases of Glimm
\cite{gli}, Dixmier \cite{dix} and Bratteli \cite{bra} and 
marked  the advent of $K$-theory in operator algebra.
For  general (nonselfadjoint) operator
algebras  of an approximately finite-dimensional nature
the situation  is more problematic and classification schemes have
usually  been restricted to those limit algebras
  $A = \indlimit A_k$ which have
intrinsic coordinates in the form of a well-defined semigroupoid.
See, for example, Poon and Wagner \cite{poo-wag}, 
Hopenwasser and Power \cite{hop-pow-1}, Muhly and Solel
\cite{muh-sol-groupoid}, and Power \cite{pow-tensor},
\cite {pow-crelle}.
In this case  the building block algebras $A_k$ are
 poset algebras, the inclusions $A_k \to A_{k+1}$ 
are regular (that is, decomposable into
multiplicity one embeddings) and the algebras $A$ 
are triangular (in the sense of Kadison and Singer \cite{kad-sin}).

In the C*-algebra direction there have been significant developments
in the last ten years in the classification of
 amenable C*-algebras using $K$-thoery invariants. This is generally
referred to as the Elliott programme; see, for example
Elliott \cite{ell-2} and D\u ad\u arlat and Eilers \cite{dar-eil}.
At the same time for nonselfadjoint operator algebras
there have been developments arising from 
 viewpoints in ring theory, representation theory
and the resolution of modules, as can be seen in
 Blecher, Muhly and Paulsen \cite{ble-muh-pau},
Muhly and Solel \cite{muh-sol} and Muhly \cite{muh}, for example.
In the present paper we generalise  the  basic 
model for C*-algebra classification
by involving representation and
embedding theory for finite-dimensional operator algebras. This leads
to   classifications of  nonselfadjoint  approximately 
finite operator algebras in terms of 
what we call  dimension module invariants.

We consider approximately finitely acting  operator algebras
as those separable algebras whose building block algebras  are 
finitely acting in the sense of having finite-dimensional generated
C*-algebras. (This gives the operator algebra category ~AFA.)
A fundamental setting occurs when these subalgebras are direct sums
of matrix algebras over a single finitely acting operator algebra
$E$. 
In particular the operator algebras determined by stationary systems
are of this form.
Here the template algebra $E$ is quite general,
notably it need not be a normed poset algebra (digraph algebra)
and the connecting homomorphisms  considered are general star-extendible
homomorphisms.
An essential point of departure with the the self-adjoint theory
is that these embeddings may not be decomposable in terms of 
multiplicity one embeddings. 

Of particular interest are those template operator 
algebras which are either 
of finite embedding type or (in some sense) of  tame embedding type.
The  former situation leads to
a classification theory
in parallel with Elliott's 
original $K_0$ classification 
of AF C*-algebras.
Away from this discrete situation
it is not the case that close embeddings are
inner unitarily equivalent but
we accommodate for this and obtain complete invariants 
by endowing the algebraic invariants with an appropriate
metric space structure.
 
The main results
 provide a framework for the 
classification of specific families of
AF operator algebras in terms of reduced dedicated invariants
and we               
give a number of applications in this direction.
For example we classify the operator algebras
which are direct limits of direct sums
of $T_r$-algebras and we consider the subcategory determined by
the $1$-decomposable (regular) 
embeddings. Also, although the usual operator 
algebra realisation
of $T_r$, the upper triangular $r $ by $r$ matrices, has
infinite embedding rank for $r \ge 3$ there are natural finitely
acting realisations
of $T_r$ of finite embedding rank, and in particular this is so for
the inflation operator algebra
$T^{max}_r$ formed by inflating over the representations
 from semi-invariant 
projections.  We determine 
the number of classes of indecomposable embeddings (the embedding rank)
from  $T^{max}_r$ into the stable algebra 
$T^{max}_r \otimes \sK $  as
\[
d(T_r^{max}) = {2r+1 \choose r+1} - (r+1).
\]
As a result we find that the classifying dimension module invariants
for limits of matrix algebras over $T_4^{max}$  (for example)
is an additive semigroup of the form $\indlimit \  \bbZ^{121}_+$ 
together with a scale and  a right  action by a 
finite multiplicative semigroup.

There are a number of other  motivations  for obtaining
nonselfadjoint generalisations of Elliott's fundamental theorem.
We show for example 
how the  abstract classification scheme  resolves
affirmatively perturbation problems of the type 
"Are close approximately finite operator algebras isomorphic ?"
For example  we prove that if limits of matricial $T_r$-algebras
(for fixed $r$)  are star-extendibly close then they are
star-extendibly isomorphic. It appears to be a 
longstanding  open question whether, in general, 
 close separable operator algebras
are isomorphic. (Examples of  Choi and Christensen \cite{cho-chr}
show that separability is necessary.)

The general topic of  perturbation  and stability for 
operator algebras, originated by Kadison and Kastler \cite{kad-kas}, 
has been well-developed for C*-algebras 
by Loring and many others. (See Loring
\cite{lor}.)
In the nonselfadjoint direction  
  a norm perturbation
theory for  reflexive operator algebras 
can be traced in  Choi and  Davidson  \cite{cho-dav},
Davidson \cite{dav},
Lance \cite{lan} and Pitts \cite{pit}.
On the other hand
the study of stability for (star-extendible) inclusions 
of nonselfadjoint
building block algebras, even in the case of digraph algebras,
is less advanced and yet highly
significant for the local structure and local 
characterisation of operator algebras.
For example 
the family of finitely acting operator algebras does not have the
perturbational stability property
of Definition 6.1 below and so it is of interest to identify 
subfamilies that do. Here we shall make use of the recent result of
Haworth \cite{haw} that the family of matricial $T_r$-algebras 
is stable in this sense.

Further motivation  for identifying invariants for limit algebras 
$A$ 
comes from the purely self-adjoint issue of the
classification of C*-subalgebra positions. 
That this connection can be made is due to the fact that
one
can often recover the pair  $\{A, A^*\}$ from the subalgebra position
\[
A \cap A^* \subseteq C^*(A)
\]
as a distinguished pair in the lattice of intermediate 
closed algebras. Thus invariants for $A$ provide invariants for the
position of $A \cap A^*$ up to C*-algebra automorphism of $C^*(A)$.
(See \cite{pow-matroid} and Section 11.)

Our starting point is the consideration  of a
finitely acting operator algebra $E$ and a family 
$\sF$ of star-extendible homomorphisms
$
\phi : E \otimes M_n \to E \otimes M_m
$
(quantum symmetries).  We associate with $\sF$  the category
$\Lim \tilde{\sF}$ of operator algebras 
of direct systems whose 
partial embeddings belong to $\sF$.
Assuming natural algebraic and analytic closure properties for 
 $\sF$,
and a certain functorial  closure property, 
we classify the algebras $A$ in ~Lim$\tilde{\sF}$ in terms of
invariants determined by $\sF$. These  invariants consist of an
ordered abelian group $
(G_\sF(A), V_\sF(A)) $
 together with a  scale
$\Sigma_{\sF}(A)$, a metric $d_A$ on the cone $V_\sF(A)$ and
the action of a metrized semiring $V_\sF$ on this cone, where $V_\sF$ 
is determined solely 
by $\sF$. In fact the additive semigroup
 $V_\sF(A)$, with its metric and 
$V_\sF$-module action, is the primary invariant
which we refer to as the
metrized dimension module of $A$. ($K_0$ is known also as the
dimension group in the case of AF C*-algebras.)
The group $G_\sF(A)$
is simply the enveloping Grothendieck group of $V_\sF(A)$.
For certain discrete families (in the sense of Definition 3.1)
  one may dispense with the metric space
structure.
In general the semiring product reflects the structure of 
compositions of embeddings between the 
building block algebras and the metric measures  the distance between 
the inner unitary equivalence classes of embeddings. 
In many examples the semiring $V_{\sF}$ is identifiable as 
 a semigroup semiring
$\bbZ_+[S]$ for
a semigroup $S$ and the $V_{\sF}$-module action reduces to an
$S$-action.

For a perturbationally stable  algebra $E$  
the invariants for the family ${\sF_E}$ of 
unrestricted star-extendible embeddings can be given more simply 
in the direct form
\[
V_{\sF_E}(A) = Hom_a(E,A \otimes \sK),~~~~
\Sigma_{\sF_E}(A) = Hom_a(E, A), ~~~~
\]
where $Hom_a(-,-)$ denotes the classes of
star-extendible homomorphisms, for approximate inner unitary
equivalence, together with the natural metric, and in this case
\[
V_{\sF_E} = Hom_a(E\otimes{\cal K}, E\otimes{\cal K})~~~~
\]
with the natural right action on $Hom_a(E,A \otimes \sK)$.
Although this manifestation  connects
with C*-algebraic  $KK$-theory, albeit with added metric,
the primary reason for our formulation of the invariants
was the identification of
$V_\sF$-action preservation as the key to the so-called existence
step in the proof of the classification theorems of Sections 4 and 5.
In earlier treatments  of partly self-adjoint operator algebas,
such as are indicated in  \cite{pow-k0},  \cite{don}, 
\cite{don-pow-2}, \cite{don-pow-3},
\cite{pow-samos}, \cite{hop-pow}
and especially \cite{pow-grot}, the existence step was resolved in 
context specific ways in terms of an identification 
of sufficient additional invariants,
such as homology and multiscales
with which to augment $K_0$. The dimension module 
formulation here however is
widely applicable.

For a proper  subfamily
$\sF \subset \sF_E$ the dimension module $V_{\sF}(A_0)$ is defined
in terms of a direct system for $A_0$ whose embeddings are morphisms
in $\sF$.
A fundamental issue which we address is whether the 
dimension module is an invariant for the algebra or, alternatively, 
depends  on the direct system in an essential way.
This is relevent to 
the classification of limits of digraph algebras with 
respect to 1-decomposable embeddings
and to the classification of standard masas up to automorphism,
including approximately inner automorphisms.
We introduce a natural property for the family $\sF$, 
which we refer to as
functoriality, and this provides a sufficient
(although not necessary) condition for  the dimension module to  be
well-defined. It is shown that 
 the 1-decomposable maps are
functorial in the case $E = T_3$  but are not functorial in the case
of complete bipartite digraph algebras. 
These results provoke the following  rather deep but natural
problems, which have connections with subfactor theory in the 
bipartite case.

Problem 1.  Determine  
those digraphs $G$  for which the 
1-decomposable digraph algebra maps 
$A(G) \otimes M_n \to A(G) \otimes M_m $, for all
$m,n$, give a functorial family. 

Problem 2. Determine the functorial
completion of nonfunctorial families and their embedding semirings.

The paper is organised as follows.
In Sections 2-5 we develop the abstract theory of classification
by metrized dimension modules. (Examples 3.5 to 3.9 indicate a
variety of finitely acting operator algebras.)
In Section 6 we discuss perturbational stability for  finitely acting
algebras and limit algebras and obtain complete invariants
for various approximately finite operator algebras.
In Sections 7 and 8 we give two applications in the case
of families $\sF_E$ of unrestricted embeddings including 
the classification of the operator algebra  limits
of inflation algebras. 
In Sections 9, 10 and 11 we discuss functorialty for constrained
embeddings and in Section 11 we indicate the
connections with  the classification of C*-subalgebra
positions (as indicated above) and with subfactors.

The main results of this paper were presented in May 1999 at the
27th Canadian Operator Theory and Operator Algebras Symposium,
in Charlottetown, PEI, and at 
the 19th  Great Plains Operator Theory Symposium 
in Ames, Iowa.

\section
{Approximately Finitely Acting Operator Algebras}

There does not appear to be an accepted 
notion of an approximately finite operator
algebra acting  on a complex Hilbert space but this has not prevented
the study of many classes of operator algebras of this nature.
We now give our preferred definitions and identify
various subcategories and supercategories. For
convenience we restrict the discussion to separable operator algebras.

Denote the category of  separable (closed)  operator algebras
by OA, with the understanding that algebras in OA are either
represented or, equivalently, are equipped with a prescribed
C*-algebra inclusion $A \to C^*(A)$.
The morphisms of OA will be taken to be the most natural
ones for subalgebras of C*-algebras per se, namely the star-extendible
homomorphisms.
These are the algebra homomorphisms that are restrictions of
C*-algebra homomorphisms between the generated C*-algebras.
In particular, if an operator
algebra $A$ is a closed union of a chain of closed subalgebras $A_k$,
then the inclusions $A_k \to A_{k+1}$ are star-extendible and
$C^*(A)$ is the closed union of the C*-algebras $C^*(A_k)$.
It is elementary but significant to note that
completely isometrically isomorphic operator algebras
need not be star-extendibly isomorphic.

Denote by OA$_n$ the subclass of algebras which can be represented,
star-extendibly, on a Hilbert space of dimension
no greater than $n$. We refer to these algebras as {\it finitely acting}
since, of course, they differ from the finite dimensional algebras
in OA. The approximately finitely acting operator
algebras are defined here 
as those which are locally approximable by
finitely acting subalgebras.
\medskip

\begin{defn}
The subcategory AFA of approximately finitely acting operator algebras
consists of those operator algebras $A$ in OA
such that for each $\epsilon > 0$ and finite family
$a_1,\dots ,a_n$ in $A$ there exists an
algebra $B$ in OA$_k$, for some $k$, and a star-extendible
embedding $\phi : B \to A$ such that
$d(a_i, \phi(B)) \le \epsilon$
for $i = 1,\dots,n.$
\end{defn}
\medskip

This definition is a local formulation which we can 
broaden further to identify an apparently wider class of algebras.
Let us write $M \subseteq_\epsilon N$ for subspaces $M, N$ of
an operator algebra if
$dist(m, N_1) \le \epsilon$ for all $ m $ in $M_1$,
where $M_1, N_1$ are the unit balls of $M, N$.

\medskip

\begin{defn}
The subcategory AFA$_\epsilon$ of almost approximately 
finitely acting operator algebras consists of 
those operator algebras $A$ in OA 
such that for each $\eta > 0$ and
$a_1,\dots ,a_n$ in $A$ there exists an
algebra $B$ in OA$_k$, for some $k$, and a C*-algebra 
embedding $\phi 
 : C^*(B) \to C^*(A)$ such that 
$\phi(B) \subseteq_\eta A$ and
$d(a_i, \phi(B)) \le \eta$
for $i = 1,\dots, n.$
\end{defn}
\medskip

On the other hand we
  write ~LIM~  for the subcategory of algebras in AFA which are
direct limit algebras $\indlimit A_k$ where
$A_k \in \  $OA$_k$ and the embeddings $A_k \to A_{k+1}$ are
star-extendible for all
$k$.

In an exactly similar way for a subfamily $\sE$ of finitely acting
operator algebras define the classes
\[
 \Lim(\sE) \subseteq \AF(\sE) \subseteq \AF(\sE)_\epsilon
\]
where   $\Lim(\sE)$ is the subclass of ~LIM~
consisting of limit algebras whose building block algebras lie in 
$\sE$. If $\sE$ is the family of elementary $E$-algebras, 
by which we mean the  finite-dimensional
operator algebras of the form
\[
E \otimes M_{n_1} \oplus \dots \oplus E \otimes M_{n_k},
\]
then we write ~AF$_E$ for 
~AF$(\sE)$.
In particular ~AF$_\mathbb{C}$ is the class of AF C*-algebras.

For the family $\sE$ of self-adjoint finite dimensional C*-algebras 
it is essentially a well-known result of Glimm that
these classes coincide \cite{gli}. The coincidence
or otherwise of analogous nonself-adjoint categories is beginning to receive
attention and  we shall make use of 
Glimm type results  in 
this direction. However we do not know if ~LIM~
$ = $~AFA~ or if ~AFA~ = ~AFA~$_\epsilon$.

More specifically  we shall be concerned with 
subcategories of ~LIM~  which derive from a given family
$\sF$ of star-extendible embeddings $\phi : A_1 \to A_2$ 
between finitely
acting operator algebras. 
It is natural to require, as we do,  the following properties of $\sF$.

(i) $\sF$ is closed under  inner unitary conjugacy.

(ii) $\sF$ is closed under compositions where such compositions
are defined.

(iii) $\sF$ is matricially stable: if $\phi \in\sF$ then
the map 
$\phi \otimes id : A_1 \otimes M_n \to A_2 \otimes M_n$   belongs to
$\sF$.

(iv) $\sF$ is sum closed : if $\phi, \psi \in\sF$ with the same 
domain algebra $A_1$ and range algebra $A_2$, then the map
$\phi \oplus \psi$ with domain $A_1$ and range $A_2 \otimes M_2$
belongs to 
$\sF$.

(v) $\sF$ is complete with respect to the norm topology.

We refer to a family $\sF$
satisfying the properties
(i) to (v)
simply as a {\it closed family}
of maps.

\medskip

\begin{defn}
Let $\sF$ be a family of star-extendible 
homomorphisms between finitely acting operator algebras
which is a closed family  in the sense above.
Then $\Lim\sF$ denotes the subset of algebras in ~LIM~ 
consisting of limit algebras of the form
$\indlimit (A_k, \phi_k)$,
where $\phi_k \in \sF$, for all $k$.
\end{defn}

\medskip

We also consider   direct sums of building block
algebras with admissible homomorphisms 
whose partial homomorphisms belong to $\sF$.
In this case we may define for a closed family $\sF$
the associated closed family $\tilde{\sF}$ of embeddings
\[
\psi : E_1 \oplus \dots \oplus E_k \to F_1 \oplus \dots \oplus F_l
\]
which have decompositions
\[
\psi  = \Sigma \oplus \psi_{ij}
\]
where each map $\psi_{ij} : E_i \to F_j$
belongs to $\sF$.
Basic examples here are the family $\sF_E$ of  all 
star-extendible embeddings
between the $E$-algebras  (algebras  
of the form $E \otimes M_n$) and the family  $\tilde{\sF}_E$
of all star-extendible embeddings
between the elementary $E$-algebras. 
More generally one may consider building block algebras that
are Morita equivalent to $E$ although we do not do so here.

\medskip

\noindent {\bf Remark 2.4}
Let $\sR$ denote the closed family of maps between digraph algebras
which are regular, or, equivalently, which are 1-decomposable.
(See Example 3.7.)  The most  natural operator algebras in LIM
are those belonging to $\Lim(\sR)$. 
In particular,
if $A$ is a closed subalgebra of an AF C*-algebra which contains 
a standard regular AF diagonal subalgebra then $A$ belongs 
to $\Lim(\sR)$.
These are the algebras which have been most extensively 
studied, particularly in the triangular case. 
(See \cite{pow-book}.)

Apart from the intrinsic interest in understanding other  operator
algebras in LIM
there are additional reasons why it is wise to admit building
block algebras which are more general than digraph operator algebras.
Firstly, limits of elementary digraph algebras with respect to
{\em non}star-extendible embeddings may yield operator algebras in LIM
which are not locally approximable by digraph algebras.
(See \cite{hef-pow} for example.)
Secondly, consider an operator algebra $A$ of the form
 \[
A = \matrrcc{C}{S}{0}{C}
\]
which acts on the Hilbert space $H \oplus H$ where $C$ is an abelian
AF C*-algebra and $S$ is a $C$-bimodule in
the algebra of compact operators, which contains no 
nonzero finite rank operators.
Although $A$ may belong to LIM, the only digraph subalgebras of $A$
are abelian and so certainly $A$ is not a star-extendible limit 
of digraph algebras.
In examples such as these it is appropriate to consider what we 
shall refer to as
inflation digraph algebras. 
These may be completely isomorphic to digraph algebras and yet not
star-extendibly so.

\section{Embedding Semirings and Embedding Rank}

We now  define the metrized embedding semiring
$V_\sF$ of a closed family of maps together with related terminology.
We also give a Krull-Schmidt type theorem for maps between finitely
acting operator algebras which  shows in particular that the
embedding semiring admits cancellation as an abelian semigroup.
At the end of the section we consider various families of finitely
acting operator algebras and their embedding semirings.

Let $E$ be a  finitely acting unital operator algebra. Write 
$Hom(E\otimes{\cal K}, E\otimes{\cal K})$
for the family of star-extendible
algebra homomorphisms
$\phi : E \otimes{\cal K} \to E\otimes{\cal K}$. We 
generally refer 
to such homomorphisms simply as maps.
Two maps $\phi, \psi$ are said to be inner equivalent or, more
emphatically, inner unitarily equivalent, if there is a unitary $u$ in the 
unitisation of $E \otimes \sK$ for which $\phi = (\Ad u) \circ \psi$.
On identifying $E$ with $E \otimes p \subseteq  E\otimes{\cal K})$,
with $p$ a rank one projection, the map $\phi$ restricts to define a
map $\phi_r : E \to E \otimes M_n$. Conversely, each such map
(a quantum symmetry in the sense of Ocneanu \cite{ocn}) determines a
unique map in $Hom(E\otimes{\cal K}, E\otimes{\cal K})$. With modest
abuse of notation we write $\phi$ for both of these maps.

Assume that  $\sF$ is  a closed family
of maps, so that $\sF$ satisfies the properties (i) - (v) of the last
section.
Write
$V_\sF$  for the set of inner unitary 
equivalence classes $[\phi]$ of the induced maps on $E \otimes {\cal
K}$ 
(with the usual convention of taking unitaries in the unitisation of
$E \otimes \sK )$.
Then $V_\sF$ is an additive
abelian group with addition given by
$[\phi] + [\psi] = [\phi \oplus \psi]$, and 
a multiplicative semigroup with product determined by composition;
 $[\phi][\psi] = [\phi \circ \psi]$.
Also, multiplication is distributive over addition and $V_\sF$ is a 
semiring with zero element and multiplicative identity.
Write  $V_E$ and also $Hom_u(E\otimes{\cal K}, E\otimes{\cal K})$ 
to indicate the
semiring determined by all maps. Then  
$V_\sF$ is a subsemiring of  $V_E$.
The enveloping ring $R_\sF$ of $V_\sF$ may be considered as the 
usual Grothendieck group
of $(V_\sF, +)$ composed of formal differences together 
with the ring product   given by
\[
([\phi]- [\psi])([\nu]-[\mu]) =
([\phi][\nu]+[\psi][\mu]) - ([\psi][\nu]+ [\phi][\mu]).
\]
In view of Theorem 3.4 below the embedding 
semiring has additive cancellation and  $V_\sF$ embeds 
injectively
in the enveloping ring.

\medskip

\begin{defn}
(i) The metrized semiring of the closed family $\sF$ is the semiring
$V_\sF$ together with the metric $d$ for which
\[
d([\psi],[\phi]) = \inf_u \|\phi - (\Ad u)\circ\psi\|
\]
where the infimum is taken over the unitary group of
the unitisation of $E \otimes \sK$.
We say that the family is discrete if this metric is discrete.

(ii) The embedding semiring $V_E$ of the operator algebra
$E$ is the metrized semiring of the family of all star-extendible 
homomorphisms
$E \otimes{\cal K} \to E \otimes{\cal K}$.
\end{defn}
\medskip

The terms in the following definitions have their
counterparts in the representaion theory 
of finite-dimensional  complex algebras.

\medskip

\begin{defn}
A (star-extendible) map between operator 
algebras is indecomposable if it cannot
be written as a direct sum of two nonzero
maps.
\end{defn}

\medskip

\begin{defn}
 A unital  operator algebra $E$ is
said to have finite embedding type
if there are only finitely many
(inner) unitary equivalence classes of indecomposable
maps in $Hom(E \otimes \sK , E \otimes \sK)$.
The embedding rank $d(E)$ is defined to be the number of these
classes.
\end{defn}

\medskip

Finite embedding type means precisely that the embedding semiring
\[
V_E = Hom_u(E \otimes \sK , E \otimes \sK)
\]
is finitely generated as an additive abelian semigroup.

Suppose now that $C^*(E)$ is simple.
 Consider  the multiplicity
of $\phi$, denoted $\mu(\phi)$, to be
the multiplicity of the
 star-extension $\tilde{\phi} : C^*(A_1) \to C^*(A_2)$.
As we have remarked in the  introduction,  it
 is a basic point of departure
with C*-algebra theory that there are indecomposable maps 
of  multiplicity greater than one.
Let us say that the  map $\phi : A_1 \to A_2$ 
is {\it $k$-decomposable} if 
it admits
a direct sum decomposition
$\phi = \phi_1 + \dots + \phi_r$ where $\mu(\phi) \le k$ for $1 \le i \le r$.
Also let us say that the finitely acting operator
algebra $E$ is $k$-decomposable if every map $E \to E \otimes \sK$
is $k$-decomposable and there are indecomposable maps of multiplicity
$k-1$.
That
$E$ may be $1$-decomposable without being of
finite embedding type in the sense below is evident in the light of 
the elementary     Toeplitz 
algebra $L_2$ of Example 3.8.

When $C^*(E)$ is simple the  semiring $V_\sF$ 
with metric ~$d$~ may be viewed as both a graded manifold,
 graded
by the multiplicity function $\mu([\phi]) = \mu(\phi)$
and as a graded semiring.
That is, there is a disjoint union
\[
V_\sF = \{0\} \sqcup V_\sF^{(1)} \sqcup  V_\sF^{(2)} \sqcup ...
\]
where the  subsets
\[
V_\sF^{(k)} = \{[\phi] : \mu([\phi]) = k\}, ~~k = 1,2, \dots ,
\]
are open-closed and satisfy
\[
V_\sF^{(k)} + V_\sF^{(l)} \subseteq V_\sF^{(k+l)},
\mbox{~\quad ~and~\quad ~} 
V_\sF^{(k)}V_\sF^{(l)} \subseteq V_\sF^{(kl)}.
\]
Also the metric $d$ is graded in the sense that if $ d([\phi], [\psi]) < 1$ 
then $[\phi], [\psi]$ belong to the same subspace $V_\sF^{(k)}.$
The fact that  $V_\sF^{(k)}$ is a locally Euclidean manifold
is exploited  in Section 6 
in the consideration of the perturbational stability
of algebras in \ Lim$\sF$.

\medskip

One can also bring into play the
2-sided module
\[
{}_FV_E = Hom_u(E \otimes \sK , F \otimes \sK)
\]
which is a right $V_E$-module and left $V_F$-module.
This is relevent in the identification of the embedding semiring $V_{E
\oplus F}$ of the  direct sum of two finitely acting operator algebras
as the matricial semiring
\[
\left[
\begin{array}{cc}
V_E & {}_EV_F\\
{}_FV_E & V_F
\end{array}
\right].
\]
However we shall not have cause here to consider such bimodules
as we restrict attention to finitely acting operator algebras
$E$ which are indecomposable. We remark that the strict counterpart to
the finite {\em representation} type of a complex 
algebra (see \cite{gab-roi} for example)
is that ${}_\bbC V_E$ have finite rank. But, as we note in the next
paragraph, this is always the case in view of the star extendibility of
the maps.

\medskip

Recall that the Krull-Schmidt
theorem for finite-dimensional complex algebra 
ensures that every finite-dimensional 
module admits a unique decomposition as a sum of finitely many
indecomposable modules.
The counterpart fact here, that every star-extendible representation
$\phi : E \to M_N $
admits, uniquely, an indecomposable decomposition
$\phi = \phi_1 \oplus \dots \oplus \phi_r $ is elementary, at 
least if $E$ is irreducible.
For in this case the star extension 
$\tilde{\phi} : C^*(E) \to M_N$, which is given a priori, 
is a map $M_r \to M_N$
which decomposes into a sum of multiplicity one embeddings.
However we wish to show, more generally, that maps between two 
finitely acting operator algebras admit
unique indecomposable decompositions, despite the fact that there may 
be uncountably many
indecomposables
and that the indecomposables need not have multiplicity one. 
Again, the proof is elementary 
C*-algebra.

\begin{thm} ({\bf Krull-Schmidt theorem})
Let $\phi : E_1 \to E_2 $
be a unital map between finitely acting operator algebras. Then  $\phi$
admits a decomposition $\phi = \phi_1 \oplus \dots \oplus \phi_r $  into
indecomposable maps and this decomposition is unique up to order and inner
unitary equivalence.
\end{thm}

\begin{proof}
Let $\{P_1, \dots , P_r \} $ be a maximal family of orthogonal projections in 
 $E_2  \cap E_2^* $ which reduces $\phi(E_1).$
Then $\phi = \phi_1 \oplus \dots \oplus \phi_r$, 
where $\phi_i(a) = P_i\phi(a)$, and this 
is an indecomposable decomposition.
Suppose that  
$\psi = \psi_1 \oplus \dots \oplus \psi_s $ is another indecomposable
decomposition, and let $Q_j = \psi_j(1).$
Then the projections 
$Q_1, \dots , Q_s  $ belong to  $E_2  \cap E_2^* $ and reduce
$\phi(E_1).$
It follows that the 
projections in the finite dimensional C*-algebra  $ C^*(P_iQ_jP_i)$
belong to  $E_2  \cap E_2^* $ and reduce  $\phi(E_1)$. Thus, 
from the maximality of the family 
 $\{P_1, \dots , P_r \} $ it follows that for each pair $i,j$ the positive operator
$P_iQ_jP_i$ is  a scalar multiple of $P_i$,
since otherwise the generated C*-algebra  contains a nonzero projection
strictly less than $P_i$.

If $\{Q_1, \dots , Q_s \} $ is also a maximal family then for each
pair $i,j$ we have  $Q_jP_iQ_j = \mu Q_j$ for some $\mu$ (depending on
the pair $i,j$). It follows there is a permutation $\pi$ such 
that  $P_i$ and $Q_{\pi(j)} $ are unitarily equivalent
in $E_2  \cap E_2^* $ and the theorem follows.

\end{proof}

There is no extensive  theory of finitely acting
operator algebras which is  ready-to-hand. Indeed, the 
problem of determining a
general classification scheme is no less involved than that of
classifying $n$-tuples of finite dimensional operators up to unitary
equivalence. 
(Some related themes are indicated in Muhly \cite{muh}.)
We can largely bypass this issue
here as our concern is directed at the role of the embedding semiring  $V_E$, and
its subsemirings, in the classification of limit algebras.  

In the remainder of this section we consider a variety of
natural template algebras $E$.  In fact all
these algebras  can be viewed as locally
partially isometric representations of the complex semigroup algebra
of a finite semigroup.
\medskip 

\noindent{\bf Example 3.5} \ \ 	
 {\bf The semiring $V_{T_2}.$}~~
 For $r> 1$ let $T_r$ denote the operator algebra
subalgebra of $M_r$ consisting of upper triangular
matrices, so that $T_r$ is spanned by the matrix units
$\{e_{i,j} : 1\le i \le j \le r\}$.
For $r = 1$, $T_r = \bbC, V_\bbC = \bbZ_+$ and $d(\bbC) = 1$.

For $r = 2$ we note that the maps
$T_2 \to T_2 \otimes M_n$ are $1$-decomposable and that there are three
classes of $1$-decomposable embeddings. Thus 
$d(T_2) = 3$ and as an additive semigroup $ V_{T_3} = \bbZ_+^3$.
Representations for the classes of indecomposables
are given by the three maps
$\theta_i : T_2 \to T_2 \otimes M_2, 0 \le i \le 2$ given by
{\small
\[
\theta_0 :
\left[\begin{array}{cc}
x&y\\
 &z
\end{array}\right]
\to
\left[
\begin{array}{cc|cc}
0&0&0&0\\
0&x&y&0\\
\hline
 & &z&0\\
 & &0&0
\end{array}\right],
\]
\[
\theta_1 :
\left[\begin{array}{cc}
x&y\\
 &z
\end{array}\right]
\to
\left[\begin{array}{cc|cc}
x&y&0&0\\
0&z&0&0\\
\hline
 & &0&0\\
 & &0&0
\end{array}\right],
~~~~~
\theta_2 :
\left[\begin{array}{cc}
x&y\\
 &z
\end{array}\right]
\to
\left[\begin{array}{cccc}
0&0&0&0\\
0&0&0&0\\
 & &x&y\\
 & &0&z
\end{array}\right],
\]
}
To see this note that for a map $\phi: T_2 \to T_2 \otimes M_n$
the partial isometry $v = \phi(e_{12})$ has $2 \times 2$ block
upper triangular form and moreover (by star-extendibility) the 
projections $v^*v$ and $vv^*$ are block diagonal. Thus
$
{\small v = \left[
\begin{array}{cc}
v_1&v_2\\
0&v_2
\end{array}
\right]. }
$
where each $v_i$ is also a partial isometry.
(Such partial isometries, whose block entries are themselves partial
isometries, are referred to as regular partial isometries with respect
to the given block structure.)
It follows that
$[\phi] = r_0[\theta_0] + r_1[\theta_1] + r_2[\theta_2]$
where $r_0 = \rank v_2$, $r_1 = \rank v_1$
and $r_2 = \rank v_3$.

In $V_{T_2}$ the three elements $[\theta_i]$ form a semigroup
$S$
and $V_{T_2}$ is the semigroup semiring $\bbZ_+[S]$.

The operator algebras of $\Lim \tilde{\sF}_{T_2}$
were classified in
 \cite{pow-k0} by augmenting the $K_0$ invariants by a partial
order on the scale of $K_0$ which derives from partial isometries
in the algebra.
Also it follows from Heffernan \cite{hef} that
\[
\Lim (\tilde{\sF}_{T_2}) = \AF(\tilde{\sF}_{T_2}) = 
\AF(\tilde{\sF}_{T_2})_\epsilon .
\]
\medskip 

\noindent{\bf Example 3.6} \ \ 	
 {\bf The semiring $V_{T_3}$.}~~
The operator algebra $T_3$ (acting on $\bbC^3$ in the usual
way)
has infinite embedding type. To see this consider
 the maps
$\phi_\alpha : T_3 \to T_3 \otimes M_3$ given by
{\small
\[
\phi_\alpha : \matrrrccc{a}{x}{z}{}{b}{y}{}{}{c}
\to
\left[
\begin{array}{ccc|ccc|ccc}
{0}&{0}&{0}&{0}&{0}&{0}&{0}&{0}&{0}\\
{0}&{0}&{0}&{0}&{0}&{0}&{0}&{0}&{0}\\
{0}&{0}&{a}&{0}&{x_1}&{0}&{x_2}&{z}&{0}\\
\hline
&&&         {a}&{x_4}&{0}&{x_3}&{0}&{z}\\
&&&         {0}&{b}&{0}&{0}&{y_1}&{y_2}\\
&&&         {0}&{0}&{0}&{0}&{0}&{0}\\
\hline
&&&&&&              {b}&{y_4}&{y_3}\\
&&&&&&              {0}&{c}&{0}\\
&&&&&&              {0}&{0}&{c}\\
\end{array}
\right]
\]
}
with
\[
\matrrcc{x_1}{x_2}{x_4}{x_3} = x\matrrcc{\alpha}{\beta}{-\beta}{\alpha},
\matrrcc{y_1}{y_2}{y_4}{y_3} = y\matrrcc{\alpha}{-\beta}{\beta}{\alpha}
\]
where $ 0 \le \beta , \alpha \le 1, \alpha^2 + \beta^2 = 1$.
If $\alpha \ne 0,1$ then $\phi_\alpha$ is an indecomposable
map with multiplicity 2.
Indeed suppose that $\phi ' + \phi '' $ is a nontrivial orthogonal
direct sum decomposition. Then $\phi'(e_{11}) = f_{33}$ or $f_{44}$
and $\phi '(e_{22}) = f_{55}$ or $f_{77}$,
where $(e_{ij})$ and $(f_{ij})$ are the underlying matrix unit
systems.
Then
\[
\|\phi'(e_{12})\| = \|\phi'(e_{11})\phi(e_{12})\phi'(e_{22})\|
= \|f_{ii}\phi'(e_{12})f_{jj}\| = |\alpha | ~~ \mbox{or} ~~ |\beta |
\]
contrary to the fact that
$\phi'$ is isometric.

The maps $\phi_\alpha$
$\phi_\gamma$
are not inner conjugate if  $\alpha \ne \gamma$ 
since, for example, the norms  of the block (1,2) entries of 
$\phi_\alpha (e_{1,2})$ and $\phi_\gamma (e_{1,2})$ differ.

On the other hand one can see that the block $(1,2)$ entry of
$\phi^{(n)}_\alpha(a)$ has norm tending to zero,
where $a \in T_3$ and  $\phi^{(n)}_\alpha$ is
the $n$-fold decomposition. From this it follows that in $V_{T_3}$ the
distance $d([\phi^{(n)}_\alpha], [\phi^{(n)}_\beta])$ tends to zero
as $n$ tendy to infinity. 
A consequence of this is that the stationary limit operator 
algebras determined
by the maps $\phi_\alpha$ for $0 < \alpha < 1$
all coincide (whilst their algebraic limits
do not).
\medskip

\noindent{\bf Example 3.7} \ \ 	
 {\bf Digraph algebras and regular maps.}~~
A finite poset $\sP = \{v_1, \dots ,v_n\}$
gives rise to a complex algebra $A(\sP)$ with $\bbC$
basis $\{e_{ij}\}$
where $i,j \in \{1,\dots,n\}$ and
$v_i \le v_j$.
A natural realisation of
$A(\sP)$ as a finitely acting operator
algebra arises when $\{e_{ij}\}$ is taken to be 
a multiplicatively closed subset  of a complete matrix unit 
system for $M_n(\bbC)$. We refer to such operator algebras as digraph
algebras. Formally a digraph algebra $A$ is a finitely 
acting operator algebra which contains a maximal abelian 
self-adjoint subalgebra (masa) of its generated C*-algebra.
The digraph $G(A)$ of $A$ is defined as the poset $\sP$, viewed as 
a transitive reflexive directed graph with no multiple edges.
The reduced digraph $G^r(A)$ is the asymmetric digraph obtained from $G(A)$
by identifying vertices in each maximal complete subgraph
$K_n \subseteq G(A)$.
In particular $G(T_2) = G(T_2(\bbC) \otimes M_n(\bbC))$,
the digraph with two vertices, two loop edges and one
proper edge.

Reciprocally, let us write $A(G)$ for the digraph operator algebra
determined by a digraph $G$ with the properties above.
Since the semirings $V_{A(G)} $ and  $V_{A(G^r)}$ coincide
we may confine attention
to reduced (antisymmetric) 
digraphs.

Let $G, H$ be connected and antisymmetric digraphs in the sense above.
Then  there are only finitely many
equivalence classes of multiplicity one
maps $\phi : A(G) \to A(H) \otimes \sK$. Indeed 
each such class corresponds to a digraph homomorphism $G \to H$.
Thus there is a natural multiplicative semigroup injection
\[
\End(G) \to Hom_u(A(G) \otimes \sK, A(G) \otimes \sK) = V_{A(G)}
\]
and a semiring injection
$ i: \bbZ_+[\End(G)] \to V_{A(G)},$ where $\End(G)$ indicates the
endomorphism semigroup of $G$.
In particular as we saw above
the map $i$ is a surjection if $A(G) = T_2$
and is not a surjection if $A(G) = T_3$.

Recall that a map $\phi : A_1 \to A_2$ between digraph algebras is 
said to be {\em regular}
if there exist masas $C_i \subseteq A_i$
such that $\phi$ maps the partial isometry 
normaliser $N_{C_1}(A_1)$ into the partial
isometry normaliser $N_{C_2}(A_2)$ \cite{pow-k0},
\cite{pow-crelle}.
The regular maps between digraph algebras are precisely those
which are 1-decomposable.
Write $\sF^{reg}_G$ for the family of regular
maps $A(G) \otimes M_n \to A(G) \otimes M_m$, for 
all $m, n$.
Then the metrized embedding
semiring for this family of maps is the semiring
$\bbZ_+[\End G]$,
with the discrete metric.
This family is closed in the sense of Section 2, as well as being
topologically closed.
\medskip 

\noindent{\bf Example 3.8} \ \ 	
{\bf Toeplitz algebras.}
Let $L_2$ be the 
finitely acting operator algebra in $M_2$ 
consisting of the matrices
\[
\left[
\begin{array}{cc}
a&b\\
0&a
\end{array}
\right].
\]
The  embedding semiring $V_{L_2}$ contains the classes of embeddings
\[
\rho_\theta \otimes i : E \otimes M_n \to E \otimes M_m
\]   
where $n < m$, $i : M_n \to M_m$ is a multiplicity one inclusion and
\[
\rho_\theta :
\left[
\begin{array}{cc}
a&b\\
0&a
\end{array}
\right]
\to
\left[
\begin{array}{cc}
a&\theta b\\
0&a
\end{array}
\right]
\]
for $|\theta | = 1$. This subset of classes, with the relative
topology,
is seen to be a homeomorph of the unit circle $S^1$.
The embeddings $\rho_\theta$ are indecomposable and so
$d(L_2) = \infty$.
Moreover it can be verified that each indecomposable map 
of multiplicity one is equivalent
to $\rho_\theta$ for some $\theta$. The only other indecomposable
is equivalent to the mulitplicity two map
$\tau : E \to E \otimes M_2$ with range in $\bbC \otimes M_2$.
It follows that $V_{L_2}$ is the abelian semiring
\[
V_{L_2} = \bbZ_+[S^1] \oplus \bbZ_+[[\tau]]
\]
with product such that
$[\rho_\theta][\tau] = [\tau][\rho_\theta] = [\tau].
$
\medskip 

\noindent{\bf Example 3.9} \ \ 	
 {\bf Quiver algebras.}
Finally let us indicate how certain finitely acting operator algebras
are derived from quivers.

Let $Q = (V, E)$ be a quiver, that is an arbitrary finite
directed graph with vertex set $V$ and edge set $E$.
A (finite directed) path $p$ of $Q$ is  either 
a trivial path $1_v$, with initial vertex and final vertex $v$,
or is a sequence $e_t\dots e_1$ of edges $e$ of $E$, for which the final
vertex of $e_{i}$ is equal to the initial vertex of 
$e_{i+1}$, for $1 \le i \le t-1$.
The path algebra $\bbC Q$ is defined to be the complex
algebra of formal linear combinations of paths with the product
defined by cocatenation of paths. This is a finite-dimensional
algebra precisely when $Q$ is acyclic.

The finite-dimensional algebras $A = \bbC Q$
with $(\rad A)^2 = 0$ correspond to quivers which 
are bipartite directed graphs. These include the poset
algebras $A(\sP)$ with this 
property, since in this case $A(\sP)$ is isomorphic to $ \bbC Q$
where $Q$ 
is the directed graph determined by $\sP$ with 
loop edges at vertices removed.
That the path algebras here are more general is due to the possibility
of multiple edges. A basic example here is the Kronecker algebra, 
which is the
complex quiver path algebra for the quiver with two edges incident
on a single vertex
and which may be realised as the subalgebra 
of $M_3$ consisting of the matrices
\[
\matrrrccc{a}{x}{y}
{0}{b}{0}{0}{0}{b}.
\]

There is a more general association of finite
dimensional complex algebras with
 quivers which is given
in terms of quotients of the path algebras of general quivers.
 Those with
$(\rad A)^2 =0$ (which have importance in the representation theory of
quivers) we may define directly as the complex algebras
$A_Q$ which are representable in terms of a matrix algebra as the set
of matrices of the form
\[
\sum_{e=(u,v)\in E} \oplus 
\bmatrix{\lambda_u}&{\lambda_e}\\{0}&{\lambda_v}  \endbmatrix
\]
where  $\lambda_u, \lambda_e, \lambda_v \in \bbC$,
$Q$ is a connected quiver and the direct sum is taken over all
directed edges $e = (u,v)$ of the quiver. 
These algebras contain the finite-dimensional path algebras 
with $(\rad A)^2 =0 $ indicated above. That they are more general is
 due to 
the admission of 
non-acyclic quivers.
In particular the $A_Q$ algebra
for the loop quiver
with a single vertex and edge
is the elementary Toeplitz algebra $L_2$.

The  presentation above gives outright a 
particular representation of 
$A_Q$ as a finitely acting operator algebra. 
Let  us denote this operator algebra  as $A_Q^{min}$. 
In particular if $Q$ is the bipartite quiver
with four edges and 4 vertices,

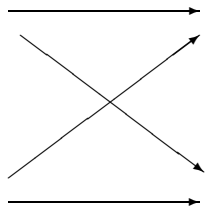
\begin{center}
\setlength{\unitlength}{3947sp}%
\begingroup\makeatletter\ifx\SetFigFont\undefined%
\gdef\SetFigFont#1#2#3#4#5{%
  \reset@font\fontsize{#1}{#2pt}%
  \fontfamily{#3}\fontseries{#4}\fontshape{#5}%
  \selectfont}%
\fi\endgroup%
\begin{picture}(1224,1224)(4789,-2773)
\thinlines
\put(4801,-1561){\vector( 1, 0){1200}}
\put(4801,-2611){\vector( 4, 3){1200}}
\put(4876,-1711){\vector( 4,-3){1152}}
\put(4801,-2761){\vector( 1, 0){1200}}
\end{picture}
\end{center}

then $A_Q^{min}$ is the operator subalgebra of $\bbC^4 \otimes M_2$
consisting of the matrices
\[
\bmatrix
{a}&{x}\\{0}&{c} 
\endbmatrix \oplus \bmatrix{a}&{y}\\{0}&{d}  \endbmatrix
\oplus \bmatrix {b}&{z}\\{0}&{d} \endbmatrix
 \oplus \bmatrix{b}&{w}\\{0}&{c} \endbmatrix .
\]
This algebra can be viewed as an inflation algebra
(in the sense of Definition 11.1) of the 4-cycle digraph algebra
$A(D_4)$. One can readily 
verify that $A_Q^{min}$ has infinite embedding
type,  in the tame sense of Example 3.8, whereas
the maximal inflation algebra $A(D_4)^{max}$, like $T_4^{max}$
has finite embedding type.

\section{Complete Invariants for Locally Finite Algebras}

We now  obtain complete
invariants for algebraic direct limit algebras whose building blocks
are $E$-algebras.
The main algebraic ideas concerning 
the existence
and uniqueness steps will reappear in the consideration
of operator algebra limits in the next section.

Let $E$ be a finitely acting operator algebra and let
$\sF$ be a family of star extendible embeddings $\phi$ between
finite dimensional
algebras of the form $E \otimes M_n$ and assume that $\sF$ 
satisfies the metric 
and algebraic closure properties
(i)-(v) of Section 2.
We do not assume that $E$ is of finite embedding rank since this
confers no particular simplification here.

Let $A \in \Lim \tilde{\sF}$ 
so that $A$ has a  presentation
\[
A = \indlimit(A_k, \phi_k)
\]
 where each $A_k$,
for $k = 1,2, \dots$, is an 
elementary $E$-algebra (that is, a finite direct sum of, let us say,
$r_k$ matrix algebras over $E$)
and where the partial embeddings 
of each map  $\phi_k$ belong to $\sF$. Also let
\[
A_0 =  \algindlimit (A_k, \phi_k)
\]
be the associated locally finite algebra.

Define $V_\sF(A_k)$ to be the right  $V_\sF$-module
which is the direct sum of $r_k$ copies of $V_\sF$
and consider  $V_\sF(A_k)$ (more functorially) as the monoid of
inner unitary equivalence classes of star extendible 
embeddings $\psi: E \to A_k \otimes {\cal K}$, where the 
partial embeddings belong to  ${\cal F}$.
Define the scale
$\Sigma_\sF(A_k)$
of $V_\sF(A_k)$ as the subset of classes $[\psi]$
for which $\psi : E \to A_k$ 
where $A_k$ is identified with $A_k \otimes \bbC p$
for some rank one projection $p$.

For each $k$ we have the induced $V_\sF$-module homomorphism
\[
\hat{\phi}_k : V_\sF(A_k) \to V_\sF(A_{k+1})
\]
given by
$\hat{\phi}_k([\psi]) = [\phi_k \circ \psi]$.
Plainly $\hat{\phi}$ respects the right $V_\sF$-action, which is to say 
that for
$[\theta]$ in $V_\sF$,
\[
\hat{\phi}_k([\psi][\theta]) = (\hat{\phi}_k([\psi]))[\theta].
\]

\begin{defn}
The dimension module of the direct system
$\{A_k, \phi_k\}$, for the family $\sF$,
is the right $V_\sF$-module
\[
V_\sF(\{A_k, \phi_k\}) = \indlimit (V_\sF(A_k), \hat{\phi}_k).
\]
\end{defn}

The direct limit here is taken in the category of additive abelian
semigroups and endowed with
 the induced right $V_\sF$-action.

Define the scale $\Sigma_\sF(A_0) $ of  $V_\sF(A_0)$ 
to be the union of the images of the scales 
$\Sigma_\sF(A_k) $ in $V_\sF(A_0)$. Writing  $G_\sF(A_0)$ for the
enveloping group of $
V_\sF(A_0)$ we obtain the scaled ordered group
\[
(G_\sF(A_0), V_\sF(A_0), \Sigma_\sF(A_0))
\]
together with the right  $V_\sF$-action on $V_\sF(A_0)$
as a tentative invariant for star-extendible isomorphism.
In  view of Proposition 3.4 the additive semigroup
$V_\sF(A_0)$ has cancellation
and so
the inclusion
$
V_\sF(A_0) \to G_\sF(A_0)$ is injective.
In fact we shall focus particularly on the possibility that the
dimension module
$ V_\sF(A_0)$ is  invariant for star-extendible homomorphisms.

The reason for caution here is that
we have no reason yet to expect functoriality
in the sense that a star extendible homomorphism
$\Phi : A_0 \to A'_0$  naturally induces an
(additive group) homomorphism  from $V_\sF(A_0)$ to   $V_\sF(A_0')$. 
Indeed we see in  the Section 11 that
 a commuting
diagram between the two direct systems which is induced by
$\Phi$  may involve morphisms which are
not associated with $\sF$. One way around this is to
require  a further  property for $\sF$, namely the 
functoriality  property of Definition 9.1. 
With a property such as this it becomes clear 
that $V_\sF(A_0)$  is indeed an invariant for the algebra, as the
notation
suggests, and is not dependent in an essential way on the particular
direct system for $A_0$.

The next theorem shows that the scaled dimension module
is a complete invariant for algebraic limit algebras determined by a
functorial family.
In the proof  we have versions of
the familiar existence and
uniqueness
steps in the construction of a commuting diagram between 
direct systems.
The existence step is relatively novel in that it 
relies on the fact that the isomorphism respects the
$V_{\cal F}$-action.
The uniqueness step is quite  elementary and closely analogous
to the self-adjoint case.

\begin{prop}
{\bf (Uniqueness.)} Let $E$ be a finitely acting operator algebra and let 
$\phi, \psi : A_1 \to A_2
$ be maps between elementary $E$-algebras with induced maps
$\hat{\phi}, \hat{\psi}$ from  $V_{\cal F}(A_1)$ to $V_{\cal F}(A_2)$.
Then $\hat{\phi}  = \hat{\psi}$ if and only if 
$\phi, \psi $ are inner  unitarily equivalent.
\end{prop}

\begin{proof}
It will be enough to establish the proposition when 
$A_1 = E \otimes M_{n_1}$ and
$A_2 = E \otimes M_{n_2}$. Let $\theta : E \otimes \sK \to 
 E \otimes \sK$ be the identity embedding.
Then 
$\hat{\phi}([\theta]) = \hat{\psi}([\theta])$
and so $[\phi \circ \theta] = [\psi \circ \theta]$ which is to say
that $[\phi] = [\psi]$ and hence that the induced maps
$\phi', \psi' : A_1 \otimes \sK  \to A_2 \otimes \sK$
are unitarily equivalent. From this it follows that
$\phi, \psi$ are unitarily equivalent.
\end{proof}

\begin{thm}
Let ${\cal F}$ be a 
family of star-extendible embeddings between operator
algebras $E \otimes M_n$, for $n=1,2,\dots,$
which is closed in the sense of Section 2, and
let $A_0, A'_0$ 
belong to $\Alglim(\tilde{\sF})$.
If ~$\Gamma$
is a $V_{\cal F}$-module isomorphism 
from $V_{\cal F}(A_0)$ to 
$V_{\cal F}(A'_0)$ then $A_0 \otimes \sK$ and $A_0' \otimes \sK$ 
are star-extendibly isomorphic.
If, moreover, $\Gamma$ gives a bijection from 
$\Sigma_\sF(A_0) $ to $\Sigma_\sF(A_0') $ then $A_0 $
and  $A_0'$ are star-extendibly isomorphic.

If 
${\cal F}$ is functorial  then the converse
of these assertions hold and the $V_{\cal F}$-module
$V_{\cal F}(-)$  is a complete invariant
for stable star-extendible isomorphism, whilst the pair
$(V_{\cal F}(-), \Sigma_\sF(-)) $ is a complete invariant
for star-extendible isomorphism.
\end{thm}

\begin{proof}
Let $A_0 = \algindlimit (A_k, \phi_k), A_0' = \algindlimit (A_k',
\phi_k')$
be the given presentations.
Consider the  $V_\sF$-module homomorphism 
$\gamma$ which is the composition
\[
 \begin{CD}
 V_\sF(A_1) @>>>  V_\sF(A_0)  \\
& \symbse{\gamma} & \symbdown{\Gamma}\\
&& V_\sF(A_0') 
\end{CD} 
\]
where $ V_\sF(A_1) \to  V_\sF(A_0)$ is the natural map.
Suppose first that $A_1 = E \otimes M_{n_1}$
and let $\theta : E \to A_1$ be the  map $a \to a \otimes p$
where  $p$ is a rank one projection.
Then $[\theta]$, the class of $\theta$ in $V_\sF(A_1)$, has image
$\gamma([\theta])$ in $V_\sF(A_0')$ which in turn coincides with the 
image of a class  $[\eta_1]$ from   $V_\sF(A_k')$, for some $k$,
under the natural map
$ V_\sF(A_k')  \to  V_\sF(A')$.
The representative  $\eta_1$ of   $[\eta_1]$
is a map $\eta_1 : E \to    A_k' \otimes M_{n_1} $, for  some $n_1$,
which in turn gives an induced map  $\eta_2$ from 
$E \otimes M_n $ to $A_k' \otimes M_{n_1n}$.

Suppose first that $A_0$ and $A'_0$ are stable algebras, so that $A_0 
= A_0 \otimes
\sK, A_0' = A_0' \otimes \sK,$ and 
$ \Sigma_\sF(A_0) =  V_\sF(A_0),
 \Sigma_\sF(A_0') =  V_\sF(A_0')$.
In particular this means that for any positive integers $k$ and $N$
one can find $l > k$ so that if
$\alpha : A_k  \to  A_l$ is the given embedding then the induced map
$\alpha \otimes \id : A_k \otimes M_N \to  A_l \otimes M_N$ is 
inner equivalent to
a map
$\gamma :  A_k \otimes M_N \to  A_l$.

Increasing $k$ if necessary it follows that we can 
replace $\eta_2$ by an inner equivalent 
 map $\eta :  E \otimes M_n \to A'_k$
such that the associated class $[\eta]$ in $ V_\sF(A_k')$
has image $\gamma([\eta])$ in $ V_\sF(A'_0)$.
It now follows that we have the factorisation
\[
 \begin{CD}
 V_\sF(A_1)&&  \\
 \symbdown{\hat{\eta}} && \symbse{\gamma}  \\
 V_\sF(A_k')  &  @>>>  & V_\sF(A_0') 
\end{CD} 
\]
Indeed, since $\gamma$ and  $\hat{\eta}$ are 
$V_{\cal F}$-module maps we have, for $[\psi]$  in 
$V_{\cal F}$,
\[
\gamma([\psi]) = \gamma([\theta][\psi]) = (\gamma([\theta])[\psi],
\]
whilst the image in $ V_\sF(A_0')$ of 
$\hat{\eta}([\psi])
$ is the image of 
\[
\hat{\eta}([\theta][\psi]) =  \hat{\eta}([\theta])[\psi] =
[\eta][\psi],
\]
and by construction, $[\eta]$ has image  $\gamma([\theta])$.
In other words, the map $\eta$ is a lifting for $\gamma$.

The case when $A_1$ has more than one summand now follows by combining
the liftings of partial embeddings.

Repeat the argument above for  the $V_\sF$-module   homomorphism
$\delta$ which is the composition
\[
 \begin{CD}
& & & V_\sF(A)  \\
&& \symbne{\delta} & \symbup{\Gamma^{-1}}  \\
 V_\sF(A_k')  &  @>>>   V_\sF(A_0') 
\end{CD} 
\]
where  $V_\sF(A_k')  \to V_\sF(A_0') $
is the natural map, to obtain a 
map $\kappa : A_k'\to A_j$
which is a lifting of $\delta$.

Since $\hat{\kappa} \circ \hat{\eta}$ is equal to the given map from
 $V_\sF(A_1)  \to V_\sF(A_j) $
it follows from Proposition 4.2 that we can replace $\delta$ by a
unitarily
equivalent map to obtain a commuting triangle. Since the process can
be repeated we obtain an infinite commuting diagram  of 
maps  between the
two given direct systems from which it follows that
$A_0$ and $A_0'$ are star-extendibly isomorphic.

Suppose now that $A_0$ and $A_0'$ are  not necessarily stable and
that $\Gamma$ preserves the scales.
Once again consider first the single summand case $A_1 = E \otimes M_{n_1}$.
If $\theta$ is as above  note  that $n_1[\theta]$ lies in 
$\Sigma_\sF(A_1) $
and so we may choose $k$ large enough so that  $n_1[\eta]$ 
is in $\Sigma_\sF(A_{k_1}'). $ It follows that there is an extension 
$\tilde{\eta}  : A_1 \to  A_{k_1}'  $ and the proof may be 
completed as before.
\end{proof}
\medskip

\section{Metrized Dimension Module Invariants}

The dimension module classification theorem for 
algebraic direct limits  
also serves to provide sufficient conditions for 
the star extendible isomorphism of operator algebra limits.
In general these conditions are   not necessary 
conditions since the closures of
$A_0$ and $A_0'$ may be isomorphic when $A_0$ and $A_0'$
are not (as we noted in Example 3.6).
In this section we obtain the appropriate invariants 
by replacing
$V_\sF(A_0)$ by its completion
under a pseudometric
induced by the metric structures on $V_\sF(A_k)$.
On the other hand, if $E$ is of finite embedding type and perturbationally stable
then we shall see in the next section that this step is unnecessary and
the functor $G : \Alglim(\sF_E) \to \Lim(\sF_E)$
is injective.

Let $E, \sF, A = \indlimit (A_k, \phi_k)$
be as given at the beginning of the last section..
Provide the $V_\sF$-modules  $V_\sF(A_k)$ with the 
natural metrics $d_k$ given by the formula of Definition 3.1,
although now $V_\sF(A_k)$ may be  a finite direct sum of 
copies of $V_\sF$.
We now define the metrized dimension module $(V_\sF(A), d)$ of the
algebra $A$ together with its presentation.

View $\hat{\phi}_k $
as a contractive metric space map from
$(V_\sF(A_k), d_k)$ to $(V_\sF(A_{k+1}), d_{k+1})$ and
define the complete metrized monoid
\[
(V_\sF(A), d) = \indlimit ((V_\sF(A_k), d_k), \hat{\phi}_k)
\]
with the limit taken in the category of metrically
complete abelian semigroups.
Explicitly this means that one forms the abelian semigroup  direct limit
\[
V_\sF^\infty(A) = \indlimit (V_\sF(A_k), \hat{\phi}_k)
\]
(which is the same as $ V_\sF(\{A_k,\phi_k\})$)  
together with the pseudometric $d$ for which
\[
d(\hat{\phi}_{k,\infty}([\psi]), \hat{\phi}_{k,\infty}([\eta]))
= \lim_{l\to \infty}
d_l(\hat{\phi}_{k,l}([\psi]), \hat{\phi}_{k,l}([\eta])).
\]
Here $\hat{\phi}_{k,\infty} : V_\sF(A_k) \to V_\sF^\infty(A)$
and $\hat{\phi}_{k, l} : V_\sF(A_k) \to V_\sF(A_l)$
are the natural homomorphisms.
Let us relax notation and write $[\eta]$     
for the typical element $\hat{\phi}_{k,\infty}([\eta])$ of
$V_\sF^\infty(A))$.
Define the equivalence relation $[\phi] \sim [\psi]$ as that for which
$d([\phi], [\psi]) = 0$.
It follows that the equivalence classes inherit a
well-defined abelian semigroup
structure. In this way  we obtain the abelian
semigroup $V_\sF^\infty(A)/\sim$, with induced metric $d$.
The completion of this metric space gives the metrized semigroup
$(V_\sF(A), d)$  which carries a unique continuous right action by $V_\sF$
which is induced from the $V_\sF$-action on $V_\sF^\infty$.
The scale $\Sigma_\sF(A)$ in $V_\sF(A)$ is defined naturally as the closure
in $V_\sF(A)$ of the natural scale
$\Sigma_{\sF}^\infty(A)/\sim$ in $V_\sF^\infty$.

\begin{thm}
Let $E$ be a finitely acting operator algebra  and let ${\cal F}$ be a 
family of star-extendible embeddings between operator
algebras $E \otimes M_n$, for $n=1,2,\dots,$
which is closed in the sense of Section 2.
Let $A, A'$ 
belong to $\Lim(\tilde{\sF})$
with direct systems determining the complete metrised semirings
$(V_{\cal F}(A), d)$ and $(V_{\cal F}(A'), d')$.
If ~$\Gamma$
is a $V_{\cal F}$-module isomorphism 
from $V_{\cal F}(A)$ to 
$V_{\cal F}(A')$ which is a bicontinuous metric space map
then $A \otimes \sK$ and $A' \otimes \sK$ 
are star-extendibly isomorphic.
If, moreover, $\Gamma$ gives a bijection from 
$\Sigma_\sF(A) $ to $\Sigma_\sF(A') $ then $A $
and  $A'$ are star-extendibly isomorphic.
\end{thm}

\begin{proof}

Let $A = \indlimit (A_k, \phi_k), A' = \indlimit (A_k',
\phi_k')$
be the given presentations.
Consider the  $V_\sF$-module homomorphism 
$\gamma$ which is the composition
\[
 \begin{CD}
 V_\sF(A_1) @>>>  V_\sF(A)  \\
& \symbse{\gamma} & \symbdown{\Gamma}\\
&& V_\sF(A') 
\end{CD} 
\]
where $ V_\sF(A_1) \to  V_\sF(A)$ is the natural map.
Suppose first that $A_1 = E \otimes M_{n_1}$
and let $\theta : E \to A_1$, be the  map $a \to a \otimes p$
where  $p$ is a rank one projection.
Then $[\theta]$, the class of $\theta$ in $V_\sF(A_1)$, has image
$\gamma([\theta])$ in $V_\sF(A_0')$.
Let $\epsilon_1 > 0.$  Choose $[\eta] $ in $V_\sF(A'_{k_1})$,
for large enough $k_1$, so that
\[
d'(\hat{\phi'}_{k_1, \infty}([\eta]),\gamma_1([\theta])) < \epsilon_1.
\]
View $\eta : E \to A'_{k_1}$ as a partially defined map
from $E \otimes M_{n}$
to $A'_{k_1}$, defined on $E \otimes \bbC p$. 
Consider first the case of stable algebras.
Then we may increase ${k_1}$ if necessary to obtain
a unique extension $\tilde{\eta}$ of $\eta$, with  $\tilde{\eta}
 : A_1 \to A'_{k_1}$.
In particular we have
\[
\hat{\tilde{\eta}}([\theta]) = [\eta],
\]
as classes in $V_\sF(A_{k_1}')$ and $\hat{\phi'}_{k_1, \infty}([\eta])
 = \hat{\phi'}_{k_1, \infty}([\tilde{\eta}])$.
By $V_\sF$-action preservation, for $[\psi]$ in $V_\sF$,
\[
\gamma_1([\psi]) = \gamma_1([\theta][\psi]) = \gamma_1([\theta])[\psi]
\]
which is $\epsilon_1$-close to \[
[\eta][\psi] = [\tilde{\eta}][\psi] = \hat{\tilde{\eta}}([\psi]),
\]
where we identify $[\eta], [\tilde{\eta}]$ with their image
classes in $V_\sF(A')$. 
This is true for all $\psi$ in $V_\sF$ and so
$\tilde{\eta}$
 is an approximate lifting
of $\gamma_1$ in the sense that the maps $\gamma_1$
and $\hat{\phi'}_{k_1,\infty} \circ \hat{\tilde{\eta}}$
are $\epsilon_1$-close as metric space maps from $V_\sF(A_1)$ to $V_\sF(A')$.
The case when $A_1$ has more than one summand now follows by combining
the liftings of partial embeddings.

Repeat the argument above for $\epsilon_2$ and the map
$\delta_1 = \Gamma^{-1} \circ \hat{\phi'}_{k, \infty}$ 
from $V_\sF(A'_{k_1})$ to 
$V_\sF(A)$
to obtain an $\epsilon_2$-approximate lifting
$\tilde{\xi}$ from
$A'_{k_1}$ to $A_{k_2}$.

Now
\[
d([\tilde{\xi} \circ \tilde{\eta}],
[\phi_{1,k_2}]) \le \epsilon_1 + \epsilon_2
\]
 and so, increasing $k_2$ we may obtain
\[ 
d_{k_2}([\tilde{\xi} \circ \tilde{\eta}],
[\phi_{1,k_2}]) \le \epsilon_1 + 2\epsilon_2,
\]
and so
 we may choose a representative
$\kappa $ in $[\tilde{\xi}]$
so that $\|\kappa \circ \tilde{\eta} - \phi_{1,k_2}\|
\le \epsilon_1 + 3\epsilon_2$.

Continue in this way, for a suitabe sequence $\epsilon_k$,
to obtain an approximately commuting diagram and the desired isomorphism.

Suppose now that $A, A'$ are not necessarily stable.
Once again consider first the single summand case $A_1 = E \otimes M_{n_1}$
and $\epsilon_1 > 0.$
Noting that $n_1[\theta]$ lies in 
$\Sigma_\sF(A_1) $
choose $k_1$ large enough so that $[\eta]$ and $n_1[\eta]$ 
are in $\Sigma_\sF(A_{k_1}') $ and
\[
d'(\hat{\phi'}_{k_1, \infty}([\eta]),\gamma_1([\theta])) < \epsilon_1.
\]
Since $n_1[\eta]$ lies in the scale there is an extension 
$\tilde{\eta}$ and the proof may be completed as before.
\end{proof}

\begin{cor}
Let $E$ be a finitely acting operator algebra and let $A, A'$
be operator algebra direct limits of $E$-algebras with respect to
star extendible embeddings. If the scaled metrized dimension
$V_E$-modules of $A$ and $A'$ are bicontinuously isomorphic 
then  $A, A'$ are star extendibly isomorphic.
\end{cor}

We note that in the case of unrestricted embeddings 
  the invariants seem to have a complexity comparable
to the limit algebras themselves.
However even the general classification 
has theoretical strength as we see 
in the next section where we obtain
perturbational stability 
for various limit algebras.
On the other hand, with prescribed embedding 
classes, or with prescribed 
building block algebras, the $V_\sF$-module action
may be replaced by a semigroup action or group
action and the  invariants become computable, as we shall see.

\section{Stability and Complete Invariants}

Let $\sE$ be a family  of finitely acting  operator algebras.
\medskip

\begin{defn}
The family $\sE$ 
has the stability property, or is perturbationally stable, if for
each algebra  $A_1$ in $\sE$ and $\epsilon >0$ 
there is a $\delta > 0$ such that to each algebra $A_2$ in $\sE$
and star-extendible
embedding
\[
\phi : A_1 \to C^*(A_2),  \ \ \ \ \  \mbox{with}\ \ \ 
\phi(A_1) \subseteq_\delta A_2
\]
there is a star-extendible embedding $\psi : A_1 \to A_2$ with 
$\|\phi - \psi\| \le \epsilon . $
\end{defn}
\medskip

In addition we say that the family is {\em uniformly
stable} if $\delta$ may be chosen independently of
the algebras in the family, and
we say that a finite-dimensional 
operator algebra $E$ is perturbationally stable if the family of
elementary $E$-algebras  is perturbationally stable.
For such an algebra $E$ it can be shown that ~AF$_E = \Lim\sF_E$.

The operator algebra $\bbC$ is well known to be 
stable and we shall see further examples  later.
For a stable algebra $E$
one has the following natural identifications for an
operator algebra inductive $A = \indlimit A_k$
of elementary $E$ algebras:
\[
Hom_a(E, A \otimes \sK) = \indlimit ~Hom_u(E, A_k \otimes \sK),\ \ 
Hom_a(E, A) = \indlimit ~Hom_u(E, A_k).
\]
Here ~$Hom_u(-,-)$ indicates the inner unitary equivalence classes of 
star-extendible homomorphisms and  ~$Hom_a(-,-)$ the 
equivalence classes of 
star-extendible homomorphisms with respect to
approximately inner automorphisms.

We can now obtain the following 
generalisation of Elliott's classification of
AF C*-algebras.

\begin{thm}
Let $E$ be a perturbationally stable  
finitely acting  operator algebra
and let $A$ 
belong to ~AF$_E$.
Then the metrized dimension module
$(V_E(A), d)$ is isometrically isomorphic to the metrized 
$V_E$-module
$Hom_a(E, A \otimes \sK)$.
Furthermore, two algebras $A, A'$ in ~AF$_E$ are star-extendibly isomorphic
if and only if there is a bicontinuous semigroup isomorphism
\[
\Gamma : Hom_a(E, A \otimes \sK) \to Hom_a(E, A' \otimes \sK),
\]
with
\[ \Gamma(Hom_a(E, A )) = Hom_a(E, A'),
\]
which respects the right action of $Hom_a(E, E \otimes \sK)$.
\end{thm}

\begin{proof}
It is clear that if$A$ and $A'$ are star-extendibly isomorphic
then there is an induced isomorphism $\Gamma$ of the invariants.
The sufficiency direction will follow from Corollary 5.2 once we show
that $V_E(A)$
is naturally isomorphic to $Hom_a(E,A\otimes\sK)$
with corresponding identification of scales.

Let $\Phi_0 : V_E^\infty(A) \to Hom_a(E, A \otimes \sK)$
be the morphism for which
$\Phi_0([\phi]) = [\phi]_{Hom}$,
where $\phi :E \to A_k \otimes M_n$
and $[\phi]$ denotes the image in $V_E^\infty(A)$ of the 
class  $[\phi]_{A_k}$ of $\phi$ in  $V_E(A_k)$,
and where $[\phi]_{Hom}$ denotes the class of $i_k \circ \phi$ where
$i_k : A_k \otimes M_n \to A \otimes \sK$ is the natural injection.
If $[\phi] = [\psi]$, then there is a sequence $u_n$ of unitaries 
in $A_{k+n}$ for which
\[
d([\phi]_{A_{k+n}}, [u_n\psi u^*_n]_{A_{k+n}}) \to 0
\]
as $n \to \infty$, which is to say that the morphisms 
$i_k  \circ \phi , i_k  \circ \psi$ from $E$ to
$A \otimes M_n$ are approximately unitarily equivalent and so 
$[\phi]_{Hom} = [\psi]_{Hom}$. For the same reasons we see 
that the metric $d(- ,-)$ on $V_E^\infty(A)$ coincides under
$\Phi_0$ with the metric on
$Hom_a(E, A \otimes \sK).$

Extending $\Phi_0$ by continuity we obtain a continuous 
isometric injection
\[
V_E(A) \to Hom_a(E, A \otimes \sK).
\]
In view of the fact that $E$ is stable this map is also surjective
and thus bijective. It also follows from the stability of $E$ that the
scale of $V_E(A)$ corresponds to $Hom_a(E, A)$.
\end{proof}

It has been shown by Haworth \cite{haw}
that the finitely acting operator algebra
$T_r$ (in $M_r$) is perturbationally stable
and hence that $\Lim\sF_{T_r} = AF_{T_r}$ for each $r \ge 1$.
Combining this with the last theorem we obtain Theorem 6.3.
The  theorem reduces to  Elliott's theorem in the case $r=1$.
For $r = 2$ star-extedible embeddings are automatically regular
(as we have seen in  Example 3.5)
and $Hom_u(T_2, T_2 \otimes \sK)$ identifies with $\bbZ^3_+$ with
discrete metric. This classification provides an alternative to
the algebraically ordered scaled ordered $K_0$ group of
Power \cite{pow-k0} and Heffernan \cite{hef}.

We refer to operator algebras in $\Lim \sF_{T_r} = \AF_{T_r}$
as approximately finite 
nest algebras of diameter $r$.

\medskip

\begin{thm}
Let $A, A'$ be approximately finite 
nest algebras of diameter $r$. Then $A, A'$
are star-extendibly isomorphic if and only if
there is a bicontinuous semigroup isomorphism
\[
\Gamma : Hom_a(T_r, A \otimes \sK) \to Hom_a(T_r, A' \otimes \sK),
\]
with
\[ \Gamma(Hom_a(T_r, A )) = Hom_a(T_r, A'),
\]
which respects the right action of $Hom_a(T_r, T_r \otimes \sK)$.
\end{thm}
\medskip

More generally the algebras of $\Lim\tilde{\sF}_{T_r}$ 
are classified in the same way.
\medskip

We now show that limit algebras determined by uniformly
stable families are themselves stable in the sense
that if $A \subseteq_\delta A'$ 
and $A' \subseteq_\delta A$
then  $A$ and $ A'$ are star extendibly isomorphic.
The hypothesis here is that the algebras are star extendibly included
in a common C*-algebra and that the Hausdorff distance between
their unit balls is no greater than $\delta$.
The key idea of the proof is to use 
the local compactness 
of $V_\sF(A)$ and $V_\sF(A')$ together with the Ascoli-Arzela
theorem to construct an isomorphism between the invariants,
and to lift this to the desired algebra
isomorphism.
\medskip

\begin{thm}
Let $E$ be a perturbationally stable
finitely acting operator algebra
and let $\sF = \sF_E$.
Let $A, A'$  be operator algebras in the class $\Lim \tilde{\sF}$
acting  on a common Hilbert space and 
suppose that $A \subseteq_\delta A'$, $A' \subseteq_\delta A$
for some $\delta > 0.$
If $\delta$ is sufficiently small then 
$A$ and $A'$ 
are star extendibly isomorphic.
\end{thm}
\medskip

\begin{proof}
Let $A = \indlimit (A_k, \phi_k),
A' = \indlimit (A_k', \phi_k')$ be presentations with 
partial embeddings 
in $\sF$. For each $k$ choose $n_k$ large enough so that
\[
A_k \subseteq_{2\delta_1} A_{n_k}'.
\]
By uniform stability we assume that $\delta_1$ is chosen so that we may
obtain in $\sF$ a star-extendible embedding
$\psi_k : A_k \to A_{n_k}'$ with
$\| \psi_k - \id\| \le 1/3$.
The map $\psi_k$ induces a map $\hat{\psi}_k$,
\[
\hat{\psi}_k : V_\sF(A_k) \to V_\sF(A'_{n_k}).
\]
Furthermore, this map is graded by multiplicity
in the sense that if $s$ is the multiplicity of $\psi_k$ then for each $r$
\[
\hat{\psi}_k^{(r)} : V_\sF(A_k)^{(r)} \to V_\sF(A'_{n_k})^{(r+s)}.
\]
The maps $\hat{\psi}_k$ are equicontinuous
metric space maps. For suppose that
$[\eta],[\nu] \in V_\sF(A_k)^{(r)}. $ Then
 \begin{align*}
d_{n_k}'(\hat{\psi}{_k([\eta])}, \hat{\psi}{_k([\nu])}) &=
\inf_{u\in U(A'_{n_k})}\|\psi_k \circ \eta - (Ad u)\circ \psi_k \circ \nu\|\\
& \le  \inf_{v\in U(A_{k})}\|\psi_k \circ (\eta - (Ad v)\circ \nu)\|\\
& \le  \inf_{v\in U(A_{k})}\|\eta - (Ad v)\circ \nu \|\\
& =  d_k([\eta],[\nu]).
 \end{align*}
 Consider the equicontinuous family of maps between the compact metric spaces
$V_\sF(A_1)^{(1)}$ and $V_\sF(A'_{n_k})^{(1+s)}$ given by the
family of restrictions
\[
\{\hat{\psi}_k|{V_\sF(A_1)}^{(1)}\}_{k=1}^\infty.
\]
By the Ascoli-Arzela theorem there is a uniformly convergent subsequence,
$\hat{\psi}_{k,1}$ say. Consider next the restrictions
\[
\hat{\psi}_{k,1}|(V_\sF(A_2)^{(1)}\cup V_\sF(A_2)^{(2)}) \to 
V_\sF(A_{n_2}')^{(1+s)} )\cup  V_\sF(A_{n_2}')^{(2+s)}
\]
and similarly obtain a uniformly convergent
subsequence $\hat{\psi}{_{k,2}}$.
Continue in this way and select a diagonal subsequence
$\hat{\psi}{_{m_k}}$ which converges uniformly on
$V_\sF(A_{j})^{(t)})$ for all $t, j$.
The limit map, $\Gamma_0$ say, inherits the properties of the maps
$\hat{\psi}{_{k}}$ in being a semiring
 homomorphism which respects the right
action
of $V_{\sF}$.
Furthermore $\Gamma_0$ is contractive and graded 
and determines a scale respecting
commuting diagram
\[
 \begin{CD}
 V_\sF(A_j) @>>>  V_\sF(A_{j+1})  \\
\symbdown{\Gamma_0}&  & \symbdown{\Gamma_0}\\
V_\sF(A_{n_j}) @>>> V_\sF(A'_{n_{j+1}})
\end{CD}
\]
where the horizontal maps are the given embeddings.
It follows that $\Gamma_0$ determines a contractive homomorphism of invariants
\[
\Gamma : (V_\sF(A), \Sigma(A),d) \to (V_\sF(A'), \Sigma(A'),d').
\]

One could appeal to the lifting arguments  of 
Theorem 4.3 at this point to construct
the desired injection but there is a direct shortcut which makes use of 
the  commuting diagrams above. By uniqueness it follows that for any
lifting
$\theta_j$ of the restriction 
\[
\Gamma_0 :   V_\sF(A_j )
\to V_\sF(A_{n_j}')
\]
there is a lifting $\theta_{j+1}$ of 
\[
\Gamma_0 :   V_\sF(A_{j+1} )
\to V_\sF(A_{n_{j+1}}')
\]
which is an extension of $\theta_j$. Thus we may 
obtain a sequence $\{\theta_j\}$ in this manner 
which determines the desired star extendible isomorphism
\end{proof}

\begin{cor}
There is a constant $\delta > 0$ such that
if $A$ and $A'$ are
approximately finite 
nest algebras of diameter $r$ and 
$A \subseteq_\delta A'$ and $A' \subseteq_\delta A$
then $A$ and $A'$ are star extendibly isomorphic.
\end{cor}

\noindent {\bf Remark 6.6}
We note that the following family of finitely acting operator algebras
does not have the stability property.

Let $\sE$ be the  family
\[
\sE = \{T^{max}_2 \otimes M_n : n=1,2,\dots \} \cup \{E_n \otimes M_m: n,m =
1,2,\dots \}
\]
where $T^{max}_2$ is the inflation algebra of matrices
\[
x = [a]\oplus \left[
\begin{array}{cc}
a&c\\
0&b
\end{array}
\right] \oplus [b]
\]
and $E_n$ is the digraph algebra
\[
E_n = 
\left[\begin{array}{cc}
D_n&S_n\\
&D_n
\end{array}
\right]
\]
where $S_n$ is the space of matrices in $M_n$ with zero
diagonal. Consider the star-extendible algebra homomorphism
\[
\phi_n : T^{max}_2 \to C^*(E_n)
\]
for which
\[
\phi_n(x) = \left[\begin{array}{cc}
aI_n&cP_n\\
&bI_n
\end{array}
\right]
\]
where $P_n$ is the rank one projection in $M_n$ all of whose entries
are $1/n$. Since $P_n$ is close to $S_n$ it is clear that
\[
\phi_n(T^{max}_2) \subseteq_{\delta_n} E_n
\]
where $\delta_n \to 0$ as $n \to \infty$. On the other hand
one can check that $\|P_n - v_n\| \ge 1/3$
for any partial isometry $v_n$ in $S_n$. (If $v_n$ were close it would 
need to be of rank one and a rank one operator with zero
diagonal is not close to $P_n$.)
It follows that
$\|\phi - \psi\| \ge 1/3$ for any star extendible
homomorphism
$\psi : T^{max}_2 \to E_n$. Thus the family $\sE$ and any family
containing $\sE$ fails to be perturbationally stable.

This example suggests that the category AFA is indeed different from
the category LIM.
\\

\noindent {\bf Remark 6.7.}\ \ 
It is plausible that the family of digraph algebras is perturbationally
stable, although not uniformly so. The following basic question also
appears to be open. Is a finitely acting operator algebra
perturbationally stable ? Choi and Davidson \cite{cho-dav} have shown
that close digraph algebras with common C*-algebra are isomorphic 
which provides some evidence for the perturabational stability of
digraph algebras. Also the 4-cycle algebras are known to give a stable
family (\cite{pow-book}) and it seems likely that the $2n$-cycle 
algebras are also perturbationally stable.

\section{Matricial $V$-algebras}

We now consider in detail the bipartite digraph algebra $E $
in $T_3$ 
 and the family $\sF \subseteq \sF_E$ of embeddings
that preserve the strictly upper triangular ideal.
This  algebra $E$ is an exceptional
complete bipartite algebra in that the maps in $\sF$ are
necessarily locally regular, although not necesarily regular. 
However the
embedding are 2-decomposable, as we see from Proposition 7.1
where we give   complete invariants for inner conjugacy. This
proposition leads to the identification of $V_\sF$ and the classification
of algebras in $\Lim (\sF)$ and $\Lim (\tilde{\sF}).$

The algebra  $E$  consists of
 matrices of the form
\[
\matrrrccc{a}{x}{y}{0}{b}{0}{0}{0}{c}
\]
and we write $V$ for  the bipartite
graph whose proper edges are indicated by the  $V$-shaped diagram

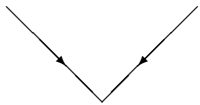
\begin{center}
\setlength{\unitlength}{3947sp}%
\begingroup\makeatletter\ifx\SetFigFont\undefined%
\gdef\SetFigFont#1#2#3#4#5{%
  \reset@font\fontsize{#1}{#2pt}%
  \fontfamily{#3}\fontseries{#4}\fontshape{#5}%
  \selectfont}%
\fi\endgroup%
\begin{picture}(1224,624)(1189,-373)
\thinlines
\put(1201,239){\vector( 1,-1){375}}
\put(2401,239){\vector(-1,-1){375}}
\put(1576,-136){\line( 1,-1){225}}
\put(1801,-361){\line( 1, 1){225}}
\end{picture}
\end{center}

Let $\varphi : E
\otimes M_{m} \rightarrow E \otimes M_{n}$ belong to $\sF$. Then 
 we may write
\[
\phi(e_{12}) = \matrrrccc{0}{v_1}{v_2}{}{0}{}{}{}{0},
~~ \phi(e_{13}) = \matrrrccc{0}{w_1}{w_2}{}{0}{}{}{}{0}.
\]
Note that since $\varphi$ is star-extendible the partial isometries
$\varphi (e_{12})$ and $\varphi (e _{13})$ have block
diagonal initial
and final projections and so $v_1, v_2, w_1, w_2$ are themselves
partial isometries.
In other words the map $\varphi$ is necessarily locally regular.

\begin{prop}
Let $\varphi, \psi : A(V) \otimes M_{m} \to  A(V) \otimes M_{n}$ 
be maps  with associated 
ordered quadruples
$\{v_1,v_2,w_1, w_2\}$, $\{x_1,x_2, y_1, y_2\}$ respectively. Then
$\varphi$ and $\psi$ are inner unitarily equivalant if and only if
\[
\rank(v_i) = \rank(x_i) , ~~ \rank(w_i) = \rank(y_i),
\]
for $i = 1,2$, and
the positive operators $w_1 w_1^*v_1 v_1^*w_1 w_1^*$, $y_1
y_1^*x_1 x_1^*y_1 y_1^*$ are unitarily equivalent.
\end{prop}

\begin{proof}
We begin the proof by putting the map $\varphi$ into a standard form.
First replace $\varphi$ by a unitary conjugate $(Adu) \circ \varphi$,
where the unitary $u$ has the form
\[
\matrrrccc{u_1}{0}{0}{}{I}{}{}{}{I},
\]
to arrange that $w_1$ and $w_2$ have standard orthogonal final
projections of the form
$$
w_1w_1^* = f_{11} +   \ldots  + f_{ss}, ~~ w_2 w_2^* =
f_{s+1,s+1} + \ldots + f_{r,r}
$$ 
where $(f_{kl})$ are the matrix units of $M_n$
and $r$ is the multiplicity of the embedding.
Next replace $\varphi$
(the new $\varphi$) by a further unitary conjugate, where the unitary
has the form
\[
\matrrrccc{I}{0}{0}{}{u_2}{}{}{}{u_3}
\]
to arrange also that the initial projections of $w_1$ and $w_2$ are
also standard orthogonal projections.
In this way arrange the
resulting standardisation of $\phi  (e_{1 3})$ to have partitioned
matrix
\[
\phi(e_{13}) = 
\left[
\begin{array}{c|cc|cc}
0&0&I_1&0&0\\
0&0&0&0&I_2\\
0&0&0&0&0\\
\hline
 &0&0& & \\
\hline
 & & &0 &0 
\end{array}
\right]
\]
where $I_1,I_2$ are sums of matrix units, $\rank(I_1) = s, 
\rank(I_2) =
t$ and $s+t \leq n$.  Since $\varphi (e_{12})$ has initial projection
  orthogonal to that of $\varphi(e_{13})$ and since $\varphi(e_{12})$
  and $\varphi(e_{13})$ have the same final projection it follows that
  $v_1$ and $v_2$ have the induced partitioned matrices
\[
v_1 = 
\left[
\begin{array}{c|c}
v_{11}&0\\
\hline
v_{12}&0\\
\hline
0&0
\end{array}
\right], ~~
 v_2 = \left[
\begin{array}{c|c}
v_{21}&0\\
\hline
v_{22}&0\\
\hline
0&0
\end{array}
\right].
\]
Note that the projections $P = v_1 v_1^*$ and $Q = w_1 w_1^*$
are block diagonal, supported in the $(1,1)$ block only, with
partitioned matrices
\[
P = v_1v^*_1 = 
\left[
\begin{array}{c|c|c}
R&S&0\\
\hline
S^*&T&0\\
\hline
0&0&0
\end{array}
\right],~~
Q = w_1w_1^* = 
\left[
\begin{array}{c|c|c}
I_1&0&0\\
\hline
0&0&0\\
\hline
0&0&0
\end{array}
\right],
\]
when $R$ is an $s \times s$ matrix and $T$ is a $t \times t$ matrix.
To simplify notation, assume without loss of generality, that $s +t =
n$, so that the last row and last column of the matrices above are not
present.
Assuming that $\psi$ satisfies the rank equality conditions in the
statement of the proposition we may also assume that $\psi(e_{13})$ is
standardised so that $\psi (e_{13}) = \varphi (e_{13})$.  We claim
that $\psi$ and $\varphi$ are inner conjugate if and only if the pair
of projections $\{P,Q\}$ in $M_n$  is unitarily equivalent
in $M_n$  to the  
pair $\{P', Q\}$, when $P' = x_1 x_1^*$. In fact the necessity
of this condition is elementary so assume that these pairs are
unitarily equivalent, by a unitary $Z$ in $M_n$.  Since $ZQZ^* = Q$ it
follows that $Z$ is block diagonal with respect to $Q$, so that
\[
Z = \matrrcc{Z_1}{0}{0}{Z_2}.
\]
Thus, conjugating $\psi$ by the unitary
\[
\matrrrccc{Z}{0}{0}{}{I_n}{}{}{}{I_n}
\]
we obtain a unitarily equivalent embedding, also in standard form,
with quadruple\\
 $\{x_1, x_2, w_1, w_2\}$ satisfying $v_1 v_1^*  = x_1
x_1^*$ (as well as $x_1^*x_1 = v_1^*v_1 , x_2^*x_2 = v_2^*v_2$).  Of
necessity $v_2v_2^* = x_2x_2^*$, since, by star extendibility
$$v_1v_1^* + v_2v_2 = \varphi(e_{12}) \varphi(e_{13})^* =
\varphi(e_{11}) = \psi(e_{11}) = \psi(e_{12})\psi(e_{12})^* = x_1 x_1^*
+ x_2 x_2^*.$$
Finally, observe that we may now conjugate $\psi$ by a unitary of the
form
\[
\left[
\begin{array}{c|cc|cc}
I_n&0&0&0&0\\
\hline
 &U_1& & & \\
 & &I_1& & \\
\hline
 & & &U_2& \\
 & & &   &I_2
\end{array}
\right]
\]
to obtain a map $\psi'$ with quadruple $\{x_1 U_1^*, x_2 U_2^*, w_1,
w_2\}$.  Thus, with the choice  $U_1 = v_1^* x_1 , U_2 = v_2^* x_2 ,
\psi' = \varphi$ and the proof is complete.
\end{proof}

Let $\varphi : A(V) \rightarrow A(V) \otimes M_n$ be as in the proof
above, with $\varphi(e_{13})$ standardised and consider again the
decomposition of the projection
\[
P = v_1v_1^* = 
\left[
\begin{array}{c | c | c}
 R & S & 0\\
 \hhline{---}
 S^* & T &0  \\
  \hhline{---}
 0 &0 &0  
\end{array}
\right]
\]
induced by $Q = w_1 w_1^*$ (with possibly 
zero rows and zero columns restored). It follows 
from the well-known spectral
picture of a pair of projections that we can further decompose the
projections $Q = w_1 w_1^*$ and $w_2 w_2^*$ so that $P$ has the form
\[
P = 
\left[
\begin{array}{c|c|c||c|c|c|c}
I_3& & & & & & \\
\hline
 &0& & & & & \\
\hline
 & &C&\sqrt{C(I-C)}& & & \\
\hline
\hline
 & &\sqrt{C(I-C)}&I-C& & & \\
\hline
 & & & &I_4& & \\
\hline
 & & & & &0& \\
\hline
&&&&&&
\end{array}
\right]
\]
where, as always,  
the unmarked entries are zero matrices and where $C$ and $I-C$
are positive invertible contractions. (See Halmos \cite{hal}.)
We also have the degenerate possibility that $C$ is  absent.

  Since $v_1
v_1^* + v_2 v_2^* = w_1 w_1^* + w_2 w_2^*$
we deduce also that $v_2 v_2^*$
has the complementary form
\[
v_2v_2^* = 
\left[
\begin{array}{c|c|c||c|c|c|c}
0& & & & & & \\
\hline
 &I_5& & & & & \\
\hline
 & &I-C&-\sqrt{C(I-C)}& & & \\
\hline
\hline
 & &-\sqrt{C(I-C)}&C& & & \\
\hline
 & & & &0& & \\
\hline
 & & & & &I_6& \\
\hline
&&&&&&
\end{array}
\right]
\]

An explicit indecomposable decomposition for $\varphi$ can now be
obtained.

Let $\theta_i, 0 \leq i \leq 3$, be the indecomposable regular
embeddings from $A(V)$ to $A(V) \otimes M_2$ whose quadruples
$\{v^{(i)}_1, v^{(i)}_2, w^{(i)}_1, w^{(i)}_2 \}$ have the rank
distributions
$$\{1,0,0,1\}, \{0,1,1,0\}, \{0,1,0,1\}, \{1,0,1,0\}.$$ 
For example we can take $\theta_2$ to be the map
\[
\theta_2 : \matrrrccc{a}{x}{y}{0}{b}{0}{0}{0}{c}
\to
\left[
\begin{array}{cccccc}
a&0&0&0&x&y\\
0&0&0&0&0&0\\
 & &0&0& & \\
 & &0&0& & \\
 & & & &b&0\\
 & & & &0&c
\end{array}
\right].
\]
A moment's reflection on the pair of 
projections $v_1 v_1^*, w_1 w_1^*$ in
the standard form of $\varphi$ above reveals  that if those projections
commute ( in which case the invertible contraction $C$ is absent) then
$\varphi$ is regular and
\[
[\varphi] = r_0[\theta_0] + r_1[\theta_1] + r_2[\theta_2] +
r_4[\theta_4]
\]
where $r_0 = \rank (I_4), r_1= \rank(I_3), r_2 = \rank(I_6), r_3 = 
\rank(I_3)$.  On the other hand if all these ranks are zero, then the
spectral diagonalisation of the invertible strict contraction $C$
gives rise to an indecomposable decomposition.  The map $\varphi$ is
indecomposable in this case if $C$ has rank one, corresponding to a
scalar $t$ in $(0,1)$.  We write  $\varphi_C$ and $\varphi_t$ 
for the maps in
these cases.  Explicitly, for a positive invertible strict contraction
$C$ of rank $n$ and for $m \geq 1$ we write $id_m \otimes \varphi_C$
for the map from
$M_m \otimes A(V)$ to $M_m \otimes A(V) \otimes M_{2n}$
given by
\[
\phi_C\left(\matrrrccc{a}{x}{y}{0}{b}{0}{0}{0}{c}\right)
=
\left[
\begin{array}{cc|cc|cc}
a\otimes I_n& &x\otimes\sqrt{C}&y\otimes I_n& x\otimes -\sqrt{I-C}
&0\\
 &a\otimes I_n & x\otimes \sqrt{I-C} &0& x\otimes\sqrt{C}&y\otimes I_n\\
\hline
 & &b\otimes I_n& & & \\
 & & &c\otimes I_n& & \\
\hline
 & & & &b\otimes I_n&  \\
 & & & & & c\otimes I_n
\end{array}
\right].
\]
This of course is the map induced by the map $\varphi_C : A(V)
\rightarrow A(V) \otimes M_{2n}$,  and the classes $[id_m \circ
\varphi_C] , [\varphi_C]$ in $V_{\cal F}$ coincide.

The next proposition is the key to exposing the semiring structure of
$V_{\cal F}$.  In view of the discussion above the assertion
follows also from the Krull Schmidt Theorem of Section 3 and the
special case of one dimensional operators $C, D$.  However the
calculation required for this scalar case is no simpler that the
calculation below.

\begin{prop}
Let $C,D$ be positive invertible strict
contractions in $M_n$ and $M_m$ respectively.  In the semiring
$V_{\cal F}$ we have $[\varphi_n][\varphi_C] = [\varphi_T]$
where $T$ is a positive operator in $M_{2nm}$  
given by
\[
T = I_2 \otimes (I_{nm} - (I_n - C) \otimes (I_n - D)).
\]
\end{prop}

\begin{proof}

  Let $\varphi_C$ have quadruple $\{v_1, v_2, w_1,
w_2 \}$ as before, let $\varphi_D$ have quadruple $\{v'_1, v'_2,
w'_2, w'_2 \}$ and consider the composition 
$(id \otimes \varphi_D) \circ \varphi_C$.  
The images of $e_{12}$ and $e_{13}$ give rise to the
quadruple $\{V_1, V_2, W_1, W_2\}$ for this composition.

Note that
$$v_1 = 1 \otimes \left[\begin{array}{ll}
                        \sqrt{C}   & 0\\
                        \sqrt{I-C} & 0\\
                        \end{array} \right], \qquad v_2 = 1 \otimes
                      \left[\begin {array}{ll}
                        -\sqrt{I-C}   & 0\\
                        \sqrt{C}      & 0\\
                        \end{array} \right],$$
$$W_1 = 1 \otimes P, \qquad w_2 = 1 \otimes P^\bot,$$
where
$$P_m  = \left[ \begin{array}{ll}
             I_m    & 0\\
             0      & 0\\
             \end{array}  \right],$$
with a similar formula for $v'_1, v'_2, w'_1, w'_2$.  Thus we may
compute
\begin{eqnarray*}
V_1  &  = & v_1 \otimes v'_1 + v_2 \otimes w'_1,\\
W_1  &  = & w_1 \otimes v'_1 + w_2 \otimes w'_1.
\end{eqnarray*}
By orthogonality
\begin{eqnarray*}
V_1 V_1^*    &  =  &  v_1v_1^* \otimes v_1' {v'_1}^* + v_2v_2^* \otimes
w'_1 {w'_1}^*,\\
            &  =  &  P_C \otimes P_D + P^\bot_C \otimes P_m,\\
W_1W_1^*     &  =  &  w_1w_1^* \otimes v'_1v_1^{'*} + w_2w_2^* \otimes
w'_1w_1^{'*},\\
            &  =  &  P_n \otimes P_D + P^\bot_n \otimes P_m.
\end{eqnarray*}
Thus $ W_1W_1^*V_1V_1^*W_1^*W_1^*$ is equal to
\[
(P_n \otimes P_D + P^\bot_n \otimes P_m)(P_C \otimes P_D + P^\bot_C \otimes P_m)
(P_n \otimes P_D + P^\bot_n \otimes P_m)
\]
\[
= (P_nP_C \otimes P_D + P^\bot_n P^\bot_C \otimes P_D P + P^\bot_n P_C \otimes P_n P_D +
   P^\bot_n P^\bot_C \otimes P_m)(P_n \otimes P_D + P^\bot_n \otimes P_m)
\]
\[
= P_n P_C P_n \otimes  P_D + P_n P^\bot_C P_n \otimes P_D P_m P_D + 
(P^\bot_n P_C P_n \otimes P_m P_D + P^\bot_n P^\bot_C P_n \otimes P_m P_D) +
\]
\[
~~(P_n P_C P^\bot_n \otimes P_D P_m +
                      P_n P^\bot_C P^\bot_n \otimes P_D P_m) +
P^\bot_n P_C P\bot_m \otimes P_m P_D P_m + P^\bot_n P^\bot_C P^\bot_n \otimes P_m.      
\]

Note that each of the bracketed terms vanishes and so we obtain
$W_1W_1^*V_1V_1^*W_1W_1^*$ as the sum of the three positive operators
{\small
\[
\left(\left[\begin{array}{cc}
C&0\\
0&0
\end{array}
(\right]
\otimes \left[\begin{array}{cc}
D&\sqrt{D(I-D)}\\
\sqrt{D(I-D)}&I-D
\end{array}
\right] +
\left[\begin{array}{cc}
I_n - C&0\\
0&0
\end{array}
\right]
\otimes \left[\begin{array}{cc}
D&\sqrt{D(I-D)}\\
\sqrt{D(I-D)}&I-D
\end{array}
\right]D\right)
\]
\[
+ \left(\left[\begin{array}{cc}
0&0\\
0&I_n -C 
\end{array}
\right] \otimes \left[\begin{array}{cc}
D&0\\
0&0
\end{array}
\right]) + 
(\left[\begin{array}{cc}
0&0\\
0&C 
\end{array}
\right] \otimes \left[\begin{array}{cc}
I_m&0\\
0&0
\end{array}
\right]\right).
\]
}

The first operator has nonzero spectral distribution
$$\{t + (1 -t)s : t \in \sigma(C), 
s \in \sigma(D)\}$$
with appropriate multiplicity whilst the sum of the second two
operators which is orthogonal to the first has nonzero spectral
distribution
$$\{(1 - t) s + t : t \in \sigma(C), 
s \in \sigma(D) \}$$
which is the same distribution.  Appealing now to Proposition 7.2
the proof is complete. 
\end{proof}

We can now identify $V_\sF$.

Define the maps $\phi_t : A(V) \to A(V) \otimes M_2$
for the values $t = 0,1$,
using the same specification as before and note that in $V_\sF$
we have
\[
[\phi_0] =  [\theta_0] + [\theta_1], ~~  [\phi_1] = [\theta_2]
+ [\theta_3]
\]
and that $\{[\phi_t]: 0 \le t \le 1\}$
is a homeomorph of the unit interval.
Let $\sF_0 \subseteq \sF$ be the (algebraically) closed
family generated by the irregular
embeddings $\phi_t $ for $0 < t < 1$,
with associated metrized semiring  ~$V_{\sF_0} ~\subseteq ~ V_\sF$.
Let $H$ be a fixed separable Hilbert space. Let $\sC_n$
be the set of strictly positive
invertible contractions of rank $n$ and let $(X_n, d)$
be the metric space where $X_n = \sC_n/\sim$ is the set
of unitary orbits of elements of $\sC_n$
and where
$d([C],[D])$ is the distance between
the unitary orbits;
\[
d([C],[D]) = inf_U \|C - UDU^*\|.
\]
More explicitly, if $\sigma (C) = \{\lambda_1, \dots ,\lambda_n\}$, 
$ \sigma (D) = \{\mu_1\dots\mu_n\}$ with repetitions reflecting
 spectral multiplicity, then the
distance
$d([C],[D])$ is the spectral distance
\[
d(\sigma (C), \sigma (D)) = inf_\pi max_i |\lambda_i - \mu_{\pi (i)}|.
\]

Set $X = \sqcup_{n=1}^\infty X_n$
with the natural induced 
metric $d$ for which $d(x_n,y_m) = 1 $ if $x_n \in X_n,
y_m \in X_m$ and $ n \ne m$.
Define a (graded) semiring structure on $X$ by defining $[C] + [D] 
= [C \oplus D]$ (using
any identification of $H \oplus H$ with $H$) and defining
$[C][D] = [C \otimes D\otimes I_2]$ (using any identification of
$H\otimes H \otimes \bbC^2$ with $H$).
Proposition 7.2 shows that the bicontinuous map $\alpha$ given by
\[
\alpha : (V_{\sF_0}, d) \to (X, d), ~~ \alpha : [\phi_C] \to [I-C]
\]
is a bicontinuous 
semiring isomorphism. The completion of the metric space $(V_{\sF_0}, d)$
is the semiring  $(V_{\sF_1}, d)$ where $\sF_1$ is the
closed family  generated by $\{[\phi_t]: t \in [0,1]\}$.

Finally note that for any class $[\psi]$ in $V_\sF$ the products
$[\psi][\theta_2],  [\psi][\theta_3], [\theta_2][\psi],
[\theta_3][\psi]  $ are 1-decomposable and so as a set
$V_\sF$ decomposes as a direct sum
\[
\bbZ_+[\theta_0] + \bbZ_+[\theta_1] + 
\bbZ_+[\theta_2] + \bbZ_+[\theta_3] + V_{\sF_0}.
\]
The semigroup structure has been determined and the topology is the
natural one consistent
with $(V_{\sF_0}, d)$.

In view of the perturbational stability of $A(V)$ and the
identification of $V_\sF$, the abstract classification theorem of
Theorem 6.2 applies and the dimension module invariants
are computable in specific cases. Given the computability of
compositions of embeddings we can obtain the following more explicit
theorem.

Let $\{C_k\}$ be a sequence of positive contractions in $M_{r_k}, $
let $ n_k = 2^kr_1r_2\dots r_k$
and let
\[
\psi_k : A(V) \otimes M_{n_k} \to A(V) \otimes M_{n_k} 
\otimes M_{2r_k}
\]
be the usual embedding with $V_{A(V)}$ class $[\psi_k] = [\phi_{C_k}]$.
Write $A(\{C_k\})$
for the unital operator algebra $\indlimit (A_k,\psi_k)$, where $A_k =
A(V) \otimes M_{n_k}$.

\begin{defn}
Two sequences $\{C_k\}, \{D_k\}$ are asymptotically
equivalent if there exist $\epsilon_k > 0$ with
$\Sigma \epsilon_k < \infty$, sequences $(n_k), (m_k)$ and
positive contractions $(X_k), (Y_k)$ such that for all $k$,
\[
d((\bigotimes _{i=n_k}^{n_{k+1}}~ (I-C_i)), ~(I-X_k) \otimes (I - Y_k)) <
\epsilon_k,
\]
\[
d((\bigotimes _{i=m_k}^{m_{k+1}} ~(I-D_i)), ~ (I-Y_k) \otimes (I - X_{k+1})) 
< \epsilon_k.
\]
\end{defn}

Here $d$ denotes the unitary orbit distance.

\begin{thm}
The operator algebras $A(\{C_k\})$ and $A(\{D_k\})$ are star extendibly
isomorphic if and only if $\{C_k\}$ and $\{D_k\}$ are aymptotically
equivalent.
\end{thm}

\begin{proof}
In view of Proposition 7.2 
if the sequences are asymptotically equivalent
then  one can construct an asymptotically commuting
diagram of unital star-extendible embeddings.
From this it follows that  $A(\{C_k\})$ and $A(\{D_k\})$ 
are star extendibly isomorphic.
On the other hand by the stability of $V$-algebras if 
$A(\{C_k\})$ and $A(\{D_k\})$ are star extendibly
isomorphic then there is an approximately commuting diagram that
implements this isomorphism. From this it follows that the sequences
are asymptotically equivalent.
\end{proof}

Write $A_C$ for the unital stationary limit algebra
$A(\{C_k\})$ where $C_k = C$ for all $k$.
It now follows readily that
if zero is not an eigenvalue of $C$ then $A_C$ is a regular limit
algebra. Also if $C, D$ are positive contractions in $M_2$
with $\sigma (C) = \{0,t\}, \sigma (D) = \{0,s\}$ with 
$0 < s,t < 1$
then $A_C$ and $A_D$ are irregular limit algebras which are
isomorphic if and only if $t = s$.

\section{Inflation algebras}

 We now classify the limit
algebras whose building block algebras are the inflation algebras
$T_r^{max} \otimes M_n$ and whose embeddings are arbitrary star
extendible embeddings.

\begin{defn}
Let $A \subseteq M_n$ be a digraph algebra and let $
\eta_i : A \to M_{n_i}, 1 \le i \le r,
$ be  unital embeddings of the form $\eta_i(a) = p_iap_i$,
where $p_i$ is a semi-invariant projection for $A$.
Then the operator algebra $A_\phi = \phi(A(G))$ in
$M_{n_1} \oplus \dots \oplus M_{n_r}$
determined by the algebra
homomorphism $\phi = \eta_1 \oplus \dots \oplus \eta_r$
is called a (regular) inflation algebra.
\end{defn}

The following example indicates why we shall be concerned here with the
case of regular inflation algebras rather than general inflation
algebras defined by contractive
representations.
Let $0 < t < 1$ and let 
$E \subseteq M_3 \oplus M_3$ be the "irregular" inflation 
algebra consisting of the matrices
\[
\matrrrccc{a}{x}{z}{}{b}{y}{}{}{c}
\oplus
\matrrrccc{a}{tx}{t^2z}{}{b}{ty}{}{}{c}.
\]
Then it can be checked that $E$ is rigid in the sense that if
$\phi: E \to E \otimes M_n$ is star extendible
and indecomposable then either $[\phi] = [id]$
or the range of $\phi$ is contained in the self-adjoint subalgebra
of $E \otimes M_n$.

Recall that the semi-invariant projections $p$ of a digraph algebra
$A$ are precisely those projections in the centre
of $A \cap A^*$ for which the correspondence $a \to pap$
determines an algebra homomorphism.
The numbers of such projections, which have the form $p_1 - p_2$
with $p_1, p_2$ invariant projections, is clearly finite.
Since repetitions of the compression embeddings
have no effect on the star-extendible isomorphism type it follows that
each digraph algebra has finitely many
nonzero regular  inflatomorphs in OA.

Write $A^{max}$ for the inflation algebra of $A$ in which all the irreducible
compression embeddings appear. Thus $T_3^{max}$
is the operator algebra of matrices of the form
\[
\matrrrccc{a}{x}{z}{0}{b}{y}{0}{0}{c}
\oplus
 \matrrcc{a}{x}{0}{b}
 \oplus \matrrcc{b}{y}{0}{c}  \oplus 
[a] \oplus [b] \oplus [c] 
\]
and its generated C*-algebra,
\[
C^*(T_3^{max}) = 
 M_3 \oplus M_2 \oplus M_2 \oplus
\bbC \oplus \bbC \oplus \bbC
\]
has maximal linear dimension amongst all the inflation algebras
of $T_3$.

We shall show that such maximal inflation algebras $T_r^{max}$ 
are of finite embedding type and compute the embedding rank
To do this we consider the  link between star-extendible embeddings
$A^{max} \to A^{max} \otimes M_N$ and nonstar-extendible embeddings
$A \to A \otimes M_N$ which are regular in the normaliser 
preservation sense
or, equivalently, which are of compression type in the
sense below.

Assume that the digraph of $A$ is connected and antisymmetric,
so that $C^*(A) = M_n$ and $A \cap A^*$ is a masa in $M_n$, which
we take to be $D_n$, the diagonal algebra.
The compression
homomorphisms $\eta_1, \dots, , \eta_r$ indicated in
Definition 8.1  also determine homomorphisms onto their ranges,
and we shall use the same notation for these maps and write 
$\eta_k : A(G) \to A(G_k)$ where $G_k$ is the 
appropriate subgraph of $G$. 

\begin{defn}
A 
{\it compression type
homomorphism}  
\[
\psi : A(G) \to A(H) \otimes M_N
\]
is an algebra
homomorphism which is unitarily equivalent to a direct
sum of {\it elementary compression type homomorphisms},
$\psi_1, \dots, \psi_s$, each of which is a composition
$\mu_k \circ \eta_k$, for some $k$,
where

\[
\begin{diagram}
\node{A(G)}   \arrow{e,t}{\eta_k}  
\node{A(G_k)} \arrow{e,t}{\mu_k} \node{A(H) \otimes M_N }
\end{diagram}
\]

\noindent and $\mu_k$ is an algebra injection arising from an 
identification
of the digraph $H_k$ with a subgraph of $G \times K_n$. Here
$K_n$ is the complete  directed graph on $N$ vertices.
\end{defn}

Up to inner unitary equivalence there are finitely
many indecomposables in the family of compression type homomorphisms
 and these
are precisely the irreducible elementary compression type
homomorphisms . When $H = G$  these indecomposables are labelled by the 
elements of the semigroup
$\Pend (G)$ of partial endomorphisms $\alpha : G \to G$ where
the domain of $G$ is a connected subgraph
determined by an interval and where $\alpha$
is a digraph homomorphism.
In Laurie and Power \cite{lau-pow} it was shown that the compression
type homomorphisms are precisely the contractive algebra
homomorphisms between digraph
algebras which are regular with respect to
some pair of masas. This has a more direct proof in the case
of $T_r$-algebras which we leave to the reader.
Using this we may obtain the following classification of maps
between the maximal inflation algebras $T^{max}_r$, $r = 1,2,\dots $.

\begin{thm}
Let $\lambda : T_r \to T_r^{max}, \kappa : T_s \to T_s^{max}$
be the canonical (nonstar-extendible)
completely isometric isomorphisms.
Let $\phi : T_r^{max} \to T_s^{max}$ be an algebra homomorphism
and let $\psi : T_r \to T_s$ be the algebra homomorphism 
$ \kappa^{-1} \circ \phi \circ \lambda$. Then the following
statements are equivalent.

(i) $\phi$ is star-extendible.

(ii) $\psi$ is of compression type.

(iii) $\psi$ is a regular contractive homomorphism.

\end{thm}

\begin{proof}
In the proof we will indicate the set of 
elementary compression type  maps for $T_r$ and $T_s$
by 
$\{\eta_i\}$
and $\{\eta_j'\}$ respectively.

Suppose first that $\psi$ is the elementary compression type  homomorphism
$\mu \circ \eta$
where
\[
\begin{diagram}
\node{T_r}\arrow{e,t}{\eta} \node{T_t}    \arrow{e,t}{\mu} 
\node{T_s}
\end{diagram}
\]
where $\eta$ is a compression type embedding determined by 
an interval projection of 
$T_r$ of rank $t$
and where $\mu $ is a multiplicity 
one star-extendible injection mapping matrix units to
matrix units.
We wish to obtain a C*-algebra extension 
\[
\tilde{\phi} : C^*(T_r^{max}) \to  C^*(T_s^{max})
\]
for the algebra homomorphism
$\phi = \kappa \circ \psi \circ \lambda^{-1}$.
Define first the restriciton
$\tilde{\phi}_{res} =
\tilde{\phi}|\eta(T_r)$, where $\eta(T_r)$ now denotes the summand
of $T_r^{max}$ corresponding to $\eta \in 
\{\eta_i\}$,  to be the map $\mu$, viewed as an algebra homomorphism
from the summand $\eta(T_r)$ to the largest
summand $\eta_{id}'(T_s) $ of $T_s^{max}$,
where $\eta_{id}' = \id$.
Since $\mu$ is star-extendible so too is this partial embedding 
$\tilde{\phi}|\eta(T_r)$, with extension given by
\[
\begin{diagram}
\node{\eta(M_r)} \arrow{e,t}{\tilde{\mu}} \node{M_s}
\end{diagram}
\]
where $\tilde{\mu}$ is the star extension of $\mu$.

We now want to fully define $\tilde{\phi}$ on all 
the other summands of C*($T_r^{max})$
so that $\tilde{\phi}$ is a C*-algebra homomorphism and 
$\tilde{\phi}(a) = 
\kappa \circ \psi \circ \lambda^{-1}$ for $a$ in 
$T_r^{max}$.
In view of the definition of $T_s^{max}$,
a matrix $b$ in $T_s^{max}$ is determined by its largest
summand, that is, by the $s \times s$ matrix summand  
$\eta_{id}'(T_s)$.
Thus the element $\kappa \circ \psi \circ \lambda^{-1}(a)$
is a direct sum of various compressions
$\eta_j'(\mu(a))$ of the $s \times s$ matrix $\mu(a)$.
The key point to note is that each such compression 
of $\mu(a)$,
which is determined by an interval of $T_s$,
can be viewed as the image
of a 
summand $\eta_{i_j}(a)$ of $a$ under a star-extendible
map.
 Let us denote this star-extension as $\tilde{\phi}_j$ ;
it is a multiplicity one C*-algebra homomorphism from 
$\eta_{i_j}(M_r)$ to $\eta_j'(M_s)$.
(If $\eta_j' = \eta'_{\id}$, the largest compression,
then $\tilde{\phi}_i = \tilde{\phi}_{res}$.)
The map $\tilde{\phi} = \sum_i \oplus \tilde{\phi}_i$
is the required extension.

We have shown  that for any elementary compression type embedding
$\psi : T_r \to T_s$ the induced algebra homomorphism
$\phi : T_r^{max} \to T_s^{max}$ is star-extendible.
It now follows that (ii) is implied by (i).

Consider now a star-extendible algebra homomorphism
$\phi :  T_r^{max} \to T_s^{max}$.
This determines a contractive  algebra homomorphism
$\psi : T_r \to T_s$ which may be viewed as the composition
\[
\begin{diagram}
\node{T_r} \arrow{e,t}{\lambda} \node{T_r^{max}}
\arrow{e,t}{\phi}  \node{T_s^{max}}
\arrow{e,t}{\pi} \node{T_s}
\end{diagram}
\]
where $\pi$ is the restriction map for the maximal dimension
summand of $C^*(T_s^{max})$.

Let $v \in T_r^{max}$ be a partial isometry. Then it is straightforward
to see that $v$ is a  unimodular sum of matrix units
and that $v = \lambda(u)$ where $u$ is a unimodular sum of matrix
units. Indeed, the maximal matrix summand of $v$ is the matrix $u$
which 
is a partial isometry along with each compression
summand.
From this it follows that $u$ is a regular partial isometry
in $T_r$.

Since $\phi$ maps partial isometries to partial isometries,
being star-extendible, it follows that $\psi(v)$ is a regular partial
isometry in $T_s$.
Thus we conclude that $\psi$ is a contractive regular homomorphism
from $T_r$ to $T_s$, and that (iii) holds.
\end{proof}

\begin{thm} 
Let $\lambda : T_r \to T_r^{max}, \kappa \otimes id : 
T_s \otimes M_n \to T_s^{max} \otimes M_n$
be the canonical (nonstar-extendible)
completely isometric isomorphisms.
Let $\phi : T_r^{max} \to T_s^{max} \otimes M_n$ 
be an algebra homomorphism
and let $\psi :  T_r \otimes M_n  \to T_s \otimes M_n$ 
be the algebra homomorphism $(\kappa \otimes id)^{-1} 
\circ \phi \circ \lambda$. Then the following
statements are equivalent.

(i) $\phi$ is star-extendible.

(ii) $\psi$ is of compression type.

(iii) $\psi$ is a regular contractive homomorphism.

\end{thm}

A proof may be given along the lines above  but for the final step.
It is not apparent
(as it is in the triangular case) that the locally regular
algebra homomorphism $\psi$ is necessarily a regular contractive
homomorphism.
That this is true for star-extendible maps was obtained recently in
Hopenwasser and Power \cite{hop-pow}
and the proof we give below is a small variation of that one.
For general digraph algebras locally regular star extendible maps need
not be regular and so the maximal triangular structure
is required in the proof.

\begin{thm}
Let $\psi : T_r \to T_s \otimes M_n$ be a contractive
algebra homomorphism.
Then the following conditions are equivalent.

(i) $\psi$ is locally regular in the sense that the image of each
matrix unit is a regular partial isometry.

(ii) $\psi$ maps regular partial isometries to regular partial
isometries.

(iii) $\psi$ is of compression type

(iv) $\psi$ is a regular contractive homomorphism.

\end{thm}

\begin{proof}
The equivalence of (i) and (ii) is elementary
and we have already noted the equivalnece of (iii)
and (iv).
Plainly (iv) implies (i).

Assume condition (i).
Let $v^{ij} =  \psi(e_{ij})$ have $s \times s$ block decomposition 
$ (v^{ij}_{pq}), $with $ 
1 \le p,q  \le s, 1 \le i,j \le r,$.
By assumption each $ v^{ij}_{pq}$ is a partial isometry. Consider a
product 
$v^{ij}v^{jk}$. The ~$(1,1)$~ block entry is
given by the sum
\[
v_{11}^{ij}v_{11}^{jk} + 
v_{12}^{ij}v_{21}^{jk} +
 \dots +
v_{1r}^{ij}v_{r1}^{jk}.
\]
Since  $v^{ij}$ is regular, 
 the
partial isometries \,$v_{11}^{ij}, \dots , v_{1r}^{ij}$ \,have orthogonal
 range
projections and so 
the operators of the sum 
have orthogonal range projections.
For similar reasons the domain projections are pairwise orthogonal.
Since, by hypothesis, the product $v^{ij}v^{jk}$ is a regular
partial isometry, it follows that the sum above is a partial isometry,
and therefore, by the orthogonality of domain and range
projections, each of
the individual products
\[
v_{11}^{ij}v_{11}^{jk}, \,v_{12}^{ij}v_{21}^{jk}, \dots , 
\,v_{1r}^{ij}v_{r1}^{jk}
\]
is a partial isometry.

Since, for example, $v_{11}^{ij}v_{11}^{jk}$ is a partial isometry
it follows that the range projection of $v_{11}^{jk}$ commutes with the
domain projection of $v_{11}^{ij}$. 
Regarding the entry operators $v^{ij}_{st}$ as identified with
operators in $M_s \otimes M_n $, 
it follows, by considering other block entries,
that for all 
$i,j,k,l,s,t,u,v$ the range projection of $v^{ij}_{st}$ commutes with the
domain projection of $v^{kl}_{uv}$. Note also that the domain
projections and the range projections commute amongst themselves.
Furthermore it is clear
that these projections commute with the 
projections in the centre of the block diagonal subalgebra of
$M_s \otimes M_n $.

Let $p_1$ be a rank one projection which is dominated by
$v_{11}^*v_{11}$. By the commutativity  there is a 
maximal family
$p_1, \dots ,p_t$  of rank
one projections satisfying
$p_iv^{i,i+1} = p_iv^{i,i+1}p_{i+1}$.
The projection $p_1 + \dots + p_t$ commutes with $\psi(T_s)$
and determines an elementary compression type embedding summand of
$\psi$. Now (iii) follows from induction.
\end{proof}

To identify $V_E$ for $E = T_r^{max}$ we need the following
combinatorial facts.

Let $[r]$ denote the totally ordered set $\{1,2, \dots ,r\}$
and let $[r,t]$ denote  the number of order preserving
functions $f : [t] \to [r]$, for
$t,r$ in $\bbN$.
Since the $[r+1,t]$ order preserving functions $g$
from $[t]$ to $[r+1]$ are partitioned into sets according
to the cardinality of $g^{-1}(1)$ we have  the recurrence identity
\[
[r+1,t] = [r,t] + [r,t-1] + \dots + [r,1] + 1.
\]
Thus
\[
[r+1,t+1] ~~=~~ [r,t+1] + [r+1,t]
\]
from which it follows that $[n,t]$ is the binomial coefficient
${n+t-1\choose t}$. In particular, if $G_r$ is the digraph for $T_r$
then we see that $\End(G_r)$ is a semigroup of
cardinality ${2r-1\choose r}$, since it is identifiable with the
semigroup of order preserving functions
$g : [r] \to [r]$.

Consider now the semigroup of partially defined order preserving functions
$h : [t] \to [r]$ whose domains are intervals.
Write $<r,t>$ for the  cardinality of this set of functions.
The set is partitioned by the cardinality of the domain of $h$, that is,
by the numbers $1, 2, \dots ,t$
and this leads to the identity
\[
<r,t>  ~~=~~ r[r,1] + (r-1)[r,2] + \dots + [r,t].
\]
Indeed, there are $r$ possible singleton domains,
$r-1$ domains of cardinality two, and so on.
Thus
\[
<r,t>  ~~=~~ r{r\choose 1} + (r-1){r+1\choose 2} +  \dots + {2r-1\choose r}.
\]
However, we have the binomial coefficient
identity
\[
{2r+1\choose r+1} = (r+1){r\choose 0} + r{r\choose 1} + 
(r-1){r+1\choose 2}
+ \dots + {2r-1\choose r},
\]
(which may be obtained by counting paths in Pascal's triangle) and so
\[
  <r,r>  ~~=~~ {2r+1\choose r+1} - (r+1).
\]

Note that the functions $h$ label the classes of indecomposable compression
type embeddings
$\eta : T_r \to T_r \otimes M_n$ (with $n \ge t$ say) and by
Theorem 8.5 these in turn correspond to the equivalence classes
of indecomposable maps $\phi : T_r^{max}  \to T_r^{max} \otimes M_n$
(for large enough $n$).

Let ${\sP}_r$ be the chain poset $\{1,\dots ,r\}$. 
Define $\Pend ({\sP}_r)$
to be  the semigroup of partially defined 
endomorphisms (monotone maps) from 
 ${\sP}_r$ to  $\sP_r$ whose domains are intervals of  $\sP_r$.
Thus the cardinality of $\Pend (\sP_r)$
is $<r,r>$ which is the embedding rank $d(T_r^{max})$. 
(It is curious that the embedding rank sequence $d_r =
d(T_r^{max})$, for $ r = 1,2, \dots $, namely,
\[
1, 7, 31, 121, 456, 1709, 6427, 24301, \dots 
\]
gives a new addition to the On-Line Encyclopedia of Integer
 Sequences, \cite{sloan}.)

\begin{thm}
Let $E = T^{max}_r$. Then the semiring $V_E$ is isomorphic to
the semiring $\bbZ_+[\Pend(\sP_r)]$
(with the discrete metric).
As an additive abelian semigroup $V_E = \bbZ_+^{d(E)}$ where
the embedding rank  is
\[
d(E) = {2r+1\choose r+1} - (r+1).
\]
\end{thm}

\begin{proof}
The first part of the theorem follows from the arguments above which show
that the indecomposable maps 
$T^{max}_r \to T^{max}_r \otimes M_n, (n > r )$ are labelled by the partially 
defined order preserving functions
$g : \{1, \dots ,r\} \to \{1, \dots ,r\}  $ whose domains are intervals.
The second assertion follows from the combinatorial discussion.
\end{proof}

The next theorem reduces the isomorphism problem for
limits of $T_r^{max}$-algebras to the structure of the embedding
semigroup. For example,  one can now compute, in principle 
(and in practice with computer aid) all the 
stationary $T_r^{max}$-algebra
limit algebras $A_{\phi}$ determined by embeddings
$\phi$ in $V_ {T_r^{max}}$ of a particular  multiplicity of low order.

\begin{thm}
Let $A$ and $A'$ be operator algebras in $\Lim(\sE)$ where $\sE$ is the
family of
$T_r^{max}$-algebras.
Let
\[
V(A) = \indlimit (\bbZ_+^{d_r}, \hat{\phi}_k), ~~
V(A') = \indlimit (\bbZ_+^{d_r}, \hat{\phi}'_k)
\]
be the dimension $V_{T_r^{max}}$-modules of $A$ and $A'$.
Then $A$ and $A'$ are stably star extendibly isomorphic
if and only if $V(A)$ and $V(A')$ are isomorphic,
and are star extendibly isomorphic 
if and only if $V(A)$ and $V(A')$ are isomorphic by a scale preserving
isomorphism.
\end{thm}

\begin{proof}
The sufficiency direction follows from Theorem 5.2 and
Theorem 8.6.
The necessity of the condition, that is, the fact that 
$V(A)$ is an invariant,
will follow from Theorem 6.2 once we show the stability of
$T_r^{max}$. However this follows readily form the 
perturbational stability of $T_r$.

Let $\alpha : 
 \to C^*(T_r^{max}
 \otimes M_n)$ be star extendible
and suppose that
\[
\alpha(T_r^{max}) ~ \subseteq_\delta ~ T_r^{max} \otimes M_n.
\]
Then consider the map $\beta : T_r \to M_r \otimes M_n$
given by $\beta = \pi \circ \alpha \circ \lambda$
where
\[
\pi : C^*(T_r^{max} \otimes M_n) \to  M_r \otimes M_n
\]
is the projection onto the largest summand.
Since $\alpha$ is an almost inclusion it follows that
\[
\beta(T_r)  ~\subseteq_\delta ~  T_r \otimes M_n.
\]
By Haworth's theorem, for $\delta$ sufficiently
small $\beta$ is close to a star extendible map
$\gamma : T_r \to T_r \otimes M_n$.
Now it follows that the map
\[
(\lambda \otimes id) \circ \gamma \circ \lambda^{-1} : 
T_r^{max} \to T_r^{max} \otimes M_n
\]
is close to $\alpha$.
\end{proof}

\section
{Functoriality and  Isoclassic Families.}

We now consider the classification problem for limit algebras 
determined by proper families of maps in $\sF_E$.

Let $\sF$ be a closed family of maps, as in Section 2.
Then $\sF$ gives rise to a number of categories
(additive $\bbC$-categories) the most elementary
of which is the category $\Sys (\sF)$ whose objects
consist of direct systems
\[
\begin{diagram}
\sA : 
\node{A_1} \arrow{e,t}{\phi_1}\node{A_2} \arrow{e,t}{\phi_2} \node{A_3}
\arrow{e} \node{\dots}
\end{diagram}
\]
with $\phi_k \in \sF$ for all $k$,
and whose morphisms are determined by commuting
diagrams with maps from $\sF$.
In particular, $\Phi : \sA \to \sA'$
is an isomorphism of $\Sys (\sF)$ if there exists a commuting diagram
of maps
\[
\begin{diagram}
  \node{A_{n_1}} \arrow[2]{e} \arrow{se}
\node[2]{A_{n_2}}
      \arrow[2]{e} \arrow{se}
\node[2]{A_{n_3}}  
 \arrow[2]{e} \arrow{se} \\
 \node[2]{A'_{m_1}} \arrow{ne}       
\arrow[2]{e}
 \node[2]{A'_{m_2}} \arrow{ne}   
\arrow[2]{e} 
  \node[2]{\makebox[1 em]{$\vphantom{A_n}$}}
\end{diagram}
\]
where the horizontal maps are compositions of the given embeddings for
$\sA , \sA'$ and the crossover maps lie  in $\sF$.

Note that the scaled dimension module $V_\sF(\sA)$ is an invariant for
morphisms in $\Sys (\sF)$ and in view of Theorem 4.3 is a complete
invariant. Thus the dimension module invariants 
resolve the isomorphism
problem for $\Sys (\sF)$.

Define the category $\Alglim ({\sF})$ whose objects 
are the operator algebras
obtained as algebraic direct limits
$A_0 = ~\algindlimit A_k$ of the systems $\sA = \{A_k, \phi_k\}$ of
~ $\Sys (\sF)$.
The morphisms of  $\Alglim (\sF)$ are the star extendible
algebra homomorphisms.
The category $\Lim (\sF)$ of closed operator algebras 
we have already indicated and
there are obvious functors
\[
\begin{diagram}
\node{Sys (\sF)} \arrow{e,t}{F} \node{\Alglim ({\sF})} 
\arrow{e,t}{G}  \node{\Lim(\sF)}
\end{diagram}
\]
However it may not be clear, even for rather elementary 
closed families, whether or not  
$F$ or $G$ induces injections (and hence
bijections) between isoclasses, that is, between the isomorphism
equivalence classes of the objects of each category.
(We leave aside here the further categorical 
questions arising from other morphisms such
as algebraic and bicontinuous morphisms. Nevertheless
see Donsig, Hudson and Katsoulis \cite{don-hud-kat} for this
consideration in the case of  regular limit algebras.)
In this connection we introduce the notion of a functorial family
of maps.
\medskip

\begin{defn}
Let $\sF$ be a closed family of maps. Then $\sF$ is said
to be functorial if for any commuting diagram
\[
\begin{diagram}
  \node{A_1} \arrow[2]{e} \arrow{se}
\node[2]{A_2}
      \arrow[2]{e} \arrow{se}
\node[2]{A_3}  
 \arrow[2]{e} \arrow{se} \\
 \node[2]{A'_1} \arrow{ne}       
\arrow[2]{e}
 \node[2]{A'_2} \arrow{ne}      
\arrow[2]{e} 
  \node[2]{\makebox[1 em]{$\vphantom{A_n}$}} 
\end{diagram}
\]
in which the horizontal maps belong to ${\sF}$
and the crossover maps are star extendible homomorphisms there
are sequences $(m_k), (n_k)
$ such that for the induced diagram
\[
\begin{diagram}
  \node{A_{n_1}} \arrow[2]{e} \arrow{se}
\node[2]{A_{n_2}}
      \arrow[2]{e} \arrow{se}
\node[2]{A_{n_3}}  
 \arrow[2]{e} \arrow{se} \\
 \node[2]{A'_{m_1}} \arrow{ne}
\arrow[2]{e}
 \node[2]{A'_{m_2}} \arrow{ne} 
\arrow[2]{e} 
  \node[2]{\makebox[1 em]{$\vphantom{A_n}$}} 
\end{diagram}
\]
the crossover maps belong to ${\sF}$.
\end{defn}
\medskip

One property which clearly leads to functoriality is the following
factorisation property.
\medskip

\begin{defn}
Let $\sF$ be a closed family of maps. Then $\sF$ is said
to have the factorisation property
if whenever $\alpha$ is a map of $\sF$ with a factorisation
$\alpha = \psi \circ \phi$ where the domains of $\phi$
and $\psi$ are  in the 
family of domain algebras for the family $\sF$, then 
$\phi$ and $\psi$ belong to $\sF$.
\end{defn}
\medskip

Plainly the functor $F$ induces an isoclass bijection if
$\sF$ is functorial.
In this case morphisms for $\Alglim (\sF)$ actually 
derive from morphisms
of $\Sys (\sF)$.  However 
the converse does not hold;
there do exist nonfunctorial families for which $F$ induces an isoclass
bijection.

The following terminology is convenient, particularly 
in the consideration 
of regular systems.
\medskip

\begin{defn}
Let $\sF$ be a closed family of maps. Then $\sF$ is said
to be an  isoclassic family if the functor $F$ is bijective.
\end{defn}
\medskip

Let $\sF^{reg}_G$ be the closed family of regular maps
$A(G) \otimes M_n \to A(G) \otimes M_m$ where $A(G)$ is 
a digraph algebra. It 
is an interesting  open question
whether $\sF^{reg}_G$ is always an isoclassic family.

For the $2n$-cycle digraph $D_{2n}$, with $n \ge 3$,
it was shown in Donsig
and Power \cite{don-pow-3} that the family of rigid regular embeddings
(those whose indecomposables derive from automorphisms
of $D_{2n}$) is a family with the factorisation
property and so functorial and isoclassic. 
Also it is shown there that the arguments admit
perturbations and that $G$ as well as $F$ gives
an isomorphism of categories.
Combining this with Theorem 4.3 we obtain the following alternative to
the 
$K_0H_1$ classification scheme of \cite{don-pow-3}.

\begin{thm}
Let $n \ge 3$ and let 
$\sA, \sA'$ be direct systems of $2n$-cycle algebras where the
embeddings belong to the family $\sF$ of maps of rigid type.
Then the following statements are equivalent

(i) $\sA, \sA'$ are isomorphic systems of $\Sys (\sF)$.

(ii) $A_0 = \algindlimit \sA , A_0' = \algindlimit \sA'$
are star extendibly isomorphic algebras.

(iii) $A = \indlimit \sA , A' = \indlimit \sA'$
are star extendibly isomorphic operator  algebras.

(iv) There is a scaled ordered group isomorphism
\[
(G_\sF(A), \Sigma_\sF(A)) \to (G_\sF(A'), \Sigma_\sF(A'))
\]
which respects the $D_{2n}$-action on the positive cones.

\end{thm}

The case $n = 2$ of 4-cycle algebras 
considered in Power \cite{pow-matroid} 
requires different methods because 
in contrast to  $ n \ge 3$ star extendible homomorphisms
need not be   locally regular.
The functor $F$ is shown  to biject isoclasses despite the
lack of the functorial property and the functor $G$ is shown to
biject isoclasses at least  in the case of odd systems.
(The even case requires a more detailed perturbational analysis.) 
In fact the 4-cycle
algebra is more naturally viewed as one of the family
of bipartite digraph algebras considered in Section 9.

It should be clear now that the two general problems
indicated in the introduction 
must be addressed in order to formulate and analyse
invariants for limit algebras.
In particular we have the specific problem of determining
 semiring $R_G$  which arises from the functorial
completion of the family of 1-decomposable embeddings of a digraph
algebra $A(G)$.

\section{Functoriality of regular $T_3$-algebra maps.}

We now show that
the family $\sF^{reg}_{T_3}$ is functorial and hence isoclassic
and we obtain a dimension module classification of the algebraic direct
limits.  
It seems quite plausible that a 
somewhat more general argument would show 
the corresponding facts
for $T_r$-algebras with $r \ge 4$.

We say that a map $\varphi : T_3 \otimes M_n \rightarrow T_3 \otimes
M_m$ is {\em $T_2$-degenerate} if $\varphi(1)$ is dominated by the sum of
two atomic interval projections of the range algebra.  Such a map is
necessarily regular for the following reasons.  Firstly the image of
each matrix unit is a regular partial isometry
in the sense that the block entries are partial isometries.
(See  the discussion of
Example 3.5.)  Secondly such locally regular maps between block upper
triangular matrix algebras are necessarily regular.  This is a special
feature of upper triangular matrix algebras given in Theorem 8.5.
Also we say that a $T_3$-algebra map $\varphi$ is of 
{\em $T_2$-character} if each indecomposable summand of $\varphi$ 
is $T_2$-degenerate.  In particular the composition
 $\varphi \circ \psi$ is also $T_2$-degenerate for an arbitrary 
map $\psi_1$ and so is regular.
Thus  ${\cal F}^{reg}_{T_3}$ is not a family with the 
factorisation property.

Consider now a direct system of $T_3$-algebras
\[
\begin{diagram}
\node{A_1} \arrow{e,t}{\phi_1}\node{A_2} \arrow{e,t}{\psi_1} \node{A_3}
\arrow{e,t}{\phi_2}\node{A_4}\arrow{e,t}{\psi_2} \node{\dots}
\end{diagram}
\]
where each map $\varphi_k$ is an irregular embedding of $T_3 \otimes
M_{n_{2k} - 1}$, and where the restriction of $\varphi_k$ to $T_2
\otimes M_{n_k}$ is regular, for each $k$.  Here $T_2$ is 
identified with the
subalgebra of $T_3$ spanned by $e_{11}, e_{12}, e_{22}$.  Assume also
that the maps $\psi_k$ are $T_2$-degenerate with $\psi_k(1)$
  contained in $T_2 \otimes M_{n_{2k+1}}$.  Then this direct system
  determines a commuting diagram isomorphism between the regular
  systems 
$
\{A_{2k - 1}, \psi_k \circ \varphi_k\}, \{A_{2k},
\varphi_{k} \circ \psi_{k - 1} \}.
$
Since the crossover maps
$\varphi_k$ are irregular this example shows that in general it is
necessary to take proper subsystems in order to establish
functoriality.  The key lemma for the proof is the following converse
to this kind of irregular factorisation.
\medskip

\begin{lma}
  Let $\varphi : A_1 \rightarrow A_2 , \psi :
A_2 \rightarrow A_3$ be maps between $T_3$-algebras and suppose that
$\psi \circ \varphi$ is regular and that $\varphi$ is irregular.
Then $\psi$ is of $T_2$-character.
\end{lma}
\medskip

\begin{proof}  We may assume that $A_1 = T_3, A_2 = T_3
\otimes M_m, A_3 = T_3 \otimes M_n$.  First note that if $\varphi$ is
irregular then for at least one of the marix units $e$ of the triple
$e_{12}, e_{23}, e_{13}$ the partial isometry $v = \varphi(e)$ has the
block form
 $$v = \left[ \begin{array}{lll}
                 a & x & z\\
 & b & y\\
                   &   & c
               \end{array} \right]$$
where the operator $b$ is not a partial isometry.  To see this 
we argue by contradiction and assume
otherwise.  Since $v^*v$ and $vv^*$ are block diagonal, the operators
$a$ and $c$ are partial isometries.  By assumption, $b$ is a partial
isometry and so it follows that the operator
$$\left[ \begin{array}{ll}
          x       & z\\
          0   & y
         \end{array} \right]$$
is a partial isometry and also has block diagonal initial and final
projections.  But this implies that $x, y$ and $z$ are partial
isometries.  We deduce then that each operator 
$ \varphi (e_{ij}), 1 \leq i
\leq j \leq 3$, is a regular partial isometry, 
which is to say  that $\varphi$
is a locally regular map.  By our remarks above this implies that
$\varphi$ is regular, contrary to hypothesis.

Since the entry $b$ of $\varphi (e)$ is not a partial isometry it
follows, by reasoning as in the last paragraph, that $x, y$ 
and $z$ are not partial isometries and in particular are nonzero operators.
Without loss of generality assume that matrix units for $A_3$
are chosen so the restriction of the map $\psi$ to the
 self-adjoint subalgebra has the form
\[
\psi : a \oplus b \oplus c \rightarrow ((a \otimes P_{11})
\oplus (b \otimes Q_{11}) \oplus (c \otimes R_{11})) \oplus
\ldots 
\oplus ((a \otimes P_{33}) \oplus (b \otimes Q_{33})
\otimes (c \otimes R_{33}))
\]
where the projections  $P   =    P_{11} + P_{22} + P_{11}, Q   =  
Q_{11} + Q_{22} + Q_{33}$  and $R   =    R_{11} + R_{22} + R_{33}$
have the same rank.
More precisely, we can remove the rows and columns of $A_3$
corresponding to the projection $(1_{A_3} - \psi(1_{A_2}))$ and obtain
the
(typical) image $\psi (v)$ in the operator matrix form
\[
\left[
\begin{array}{ccc|ccc|ccc}
 a \otimes P_{11} & x \otimes X_{11}&z \otimes Z_{11}
 &0&x\otimes X_{12}&z\otimes Z_{12}&0&x\otimes X_{13}&z\otimes Z_{13}  \\
 0 &b \otimes Q_{11}&y \otimes
 Y_{11}&0&0&y\otimes Y_{12}&0&0&y\otimes Y_{13}  \\
 0 &0&c \otimes R_{11}&0&0&0&0&0&0  \\
 \hhline{---------}
   & & &a \otimes P_{22} & x \otimes X_{22}&z \otimes
Z_{22}&0&x\otimes X_{23}&z\otimes Z_{23}  \\
   & & &0&b\otimes Q_{22}&y\otimes Y_{22}&0&0&y\otimes Y_{23}  \\
   & & &0&0&c\otimes R_{22} &0&0&0  \\
 \hhline{---------}
   & & & & & & a \otimes P_{33} & x \otimes X_{33}&z \otimes Z_{33}  \\  
   & & & & & & 0&b\otimes Q_{33}&y\otimes Y_{33}  \\
   & & & & & & 0&0&c\otimes R_{33} 
\end{array}
\right].
\]

We shall now show that the map $\psi$ is locally regular.
Since the composition $\psi \circ \varphi$ is assumed to be regular
the matrix above, for the element 
$v = \varphi(e)$, is a regular partial isometry.  In particular the
$(2,2)$ block entry and the $(2,3)$ block entry have orthogonal ranges
and so
\[
\left[
\begin{array}{ccc}
0&0&0\\ x^*\otimes X_{23}^*&0&0\\z^* \otimes Z_{23}^* & y^*
\otimes Y_{23}^*&0
\end{array}\right]
\left[
\begin{array}{ccc}
a \otimes P_{22} & x \otimes X_{22}&z \otimes Z_{22}\\
0&b\otimes Q_{22}&y\otimes Y_{22}\\
0&0&c\otimes R_{22}
\end{array}\right] = 0.
\]
Thus $x^*x \otimes {X_{23}}^*X_{22} = 0$, and so ${X_{23}}^*X_{22} = 0$.
If $p$ is a rank one projection in $M_m$ then the partial isometry
$\psi (e_{12} \otimes p)$ has block diagonal initial and final
projections and, after removal of block rows and columns of zeros, has
the $3 \times  3$ block matrix form $(p \otimes X_{ij})$, where $X_{ij} = 0$
for $i > j$.  
Since ${X_{23}}^*X_{22} = 0$  it follows from 
the block diagonality of the range projections
 that
$X_{22}$   and $X_{23}$ are partial isometries.
Reasoning as before it follows that
 $\psi ( e_{12} \otimes p)$ is a regular partial
isometry.

Similarly, since the $(2,3)$ block and the $(3,3)$ 
block of $\psi (v)$
 have orthogonal initial projections it follows that
\[
\left[
\begin{array}{ccc}
0&x\otimes X_{23}&z\otimes Z_{23}\\
0&0&y\otimes Y_{23}\\
0&0&0
\end{array}
\right]
\left[
\begin{array}{ccc}
a^*\otimes P_{33}&0&0\\ x^*\otimes X_{33}^*&b^*\otimes Q_{33}&0\\
z^* \otimes Z_{33}^* & y^* \otimes Y_{33}^*& c^*\otimes R_{33}
\end{array}\right] = 0.
 \]
In particular $yy^* \otimes Y_{23} Y_{33}^* = 0$ and so it follows,
as before, that $\psi (e_{23} \otimes p)$ is a regular partial
isometry.

If $U$ and $V$ are regular partial isometries with $U^*U = VV^*$ then it
need not be the case that the partial isometry $UV$ is regular.  Thus
we need additional argument in order to see that $\psi (e_{13} \otimes
p)$ is a regular partial isometry.
Returning once more to the regularity of $\psi (v)$ 
and the orthogonality of the range projections of the (2,2) and the
(2,3) block entries, we see that
\[
z^*z \otimes Z^*_{23}Z_{23} + y^*y \otimes Y^*_{23}Y_{23} = 0.
\]
Since $Y^*_{23}Y_{23} = 0$, by the regularity of the partial isometry
$\psi(e_{23})$ it follows, since $z \ne 0$, that 
$Z^*_{23}Z_{23} = 0$. Since $\psi(e_{13})$ has block diagonal final
projection this implies that $Z_{22}$ and $Z_{23}$ are partial
isometries, and so, as before, $\psi(e_{13})$ is a regular partial
isometry. Thus, by our earlier remarks, $\psi$ is regular.

Suppose now that $\psi_1$ is an indecomposable summand of $\psi$
which is necessarily of multiplicity one.  By
the Krull-Schmidt theorem, Theorem 3.4, 
indecomposable decompositions are unique and  from this it follows
that since
$\psi \circ \varphi$ is regular so too is $\psi_1 \circ \varphi$.  If
$\psi_1$ is not $T_2$-degenerate then $\psi (v)$ is not a regular
partial isometry, contrary to the regularity of $\psi_1 \circ
\varphi$.  It follows then that $\psi$ is of $T_2$-character.

\end{proof}

\begin{thm}
Let $\sF^{reg}_{T_3}$ be the family of regular embeddings
between $T_3$-algebras. Then $\sF^{reg}_{T_3}$
and $\tilde{\sF}^{reg}_{T_3}$ are isoclassic families.
\end{thm}

\begin{proof}
We first note the immediate consequence of Lemma 7.4
that if $\sA = \{A_k, \alpha_k\}$ and $\sA' = \{A_k', \beta_k\}$
are $T_3$-algebra systems for which none of the embeddings
$\alpha_k$, $\beta_k$ and their system compositions
 are of $T_2$-character then a commuting diagram
isomorphism between $\sA$ and $\sA'$ is necessarily regular.

In general consider the commuting diagram 
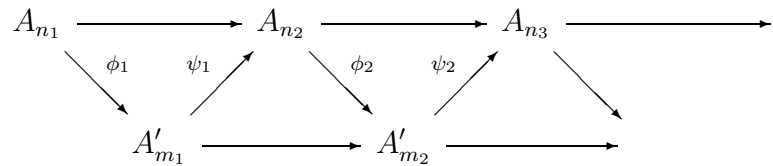
\[
\begin{diagram}
  \node{A_{n_1}} \arrow[2]{e} \arrow{se,t}{\phi_1} \node[2]{A_{n_2}}
      \arrow[2]{e} \arrow{se,t}{\phi_2} \node[2]{A_{n_3}}  
 \arrow[2]{e} \arrow{se} \\
 \node[2]{A'_{m_1}} \arrow{ne,t}{\psi_1} \arrow[2]{e}
 \node[2]{A'_{m_2}} \arrow{ne,t}{\psi_2} \arrow[2]{e} 
  \node[2]{\makebox[1 em]{$\vphantom{A_n}$}}
\end{diagram}
\]
with $\phi_k, \psi_k$ star extendible for all $k$.
Suppose moreover that infinitely
many of the maps $\phi_k$ are irregular.
Replacing the systems by subsystems we may assume that all these
maps are irregular. By the
lemma all the maps $\psi_k$ are of $T_2$-character.
Since $\phi_k \circ \psi_{k-1}$ is regular it must be 
that the range of  the regular map 
$\psi_{k-1}$ does not meet those off-diagonal blocks which
contain rank one matrix units $e$ for which $\phi_k(e)$ is an irregular
partial isometry.
This implies that 
$\phi_k \circ \psi_{k-1} \circ \phi_{k-1}$ is locally regular.
Indeed the range of $\psi_{k-1}$ meets more blocks of $A_k$ than does
the range of 
(the regular map) $\psi_{k-1} \circ \phi_k$.
Since the triple composition is locally regular it is regular and it
now follows that there is a commuting subdiagram of regular maps.
Thus $\sF^{reg}_{T_3}$ is functorial and hence an isoclassic family.
The argument is the same for $\tilde{\sF}^{reg}_{T_3}$. 
\end{proof}

We can now  see that the locally finite algebras determined by regular
embeddings of $T_3$-algebras have well-defined
dimension module
invariants and moreover are classified by these
invariants. The classification of the corresponding
operator algebras requires an peturbational version of the last lemma.

\begin{cor}
Let $A_0, A'_0$
belong to $\Alglim ({\sF}^{reg}_{T_3})$
and let
\[
V(A_0) = \indlimit (\bbZ_+^{10}, \hat{\phi}_k),~~ 
V(A_0') = \indlimit (\bbZ_+^{10}, \hat{\phi}'_k)
\]
be their dimension modules with right action from the semigroup
$V_{{\sF}^{reg}_{T_3}} = \bbZ^{10}_+$ determined by 
their defining
 direct systems. 
Then $A_0$ and $A_0'$ are stably star extendibly isomorphic
if and only if the dimension modules
 $V(A_0)$ and $V(A_0')$ are isomorphic,
and are star extendibly isomorphic 
if and only if $V(A_0)$ and $V(A_0')$ are isomorphic by a 
scale preserving
isomorphism.
\end{cor}

\section{Bipartite Digraphs and Nonfunctoriality.}

We now show the nonfunctoriality of regular embeddings of complete
bipartite digraph algebras. This is done by constructing
 pairs of  high multiplicity
indecomposable embeddings
with compositions that are 1-decomposable.
We also indicate connections with subfactors and positions
of  self-adjoint subalgebras of C*-algebras.

Let $G_{n,m}$ be the complete bipartite digraph whose digraph algebra is
\[
A(G_{n,m}) = \left[\begin{array}{cc}
\bbC^n&M_{n,m}\\
0&\bbC^m
\end{array}\right]
\]
where $M_{n,m}$ is the ${\bbC^n}-{\bbC^m}$-bimodule of $n \times n$ complex
matrices. Also write $G_n$ for $G_{n,n}$. In particular $A(G_1) = T_2,
A(G_2)$ is the 4-cycle algebra and $A(G_{1,2})$ is the $V$-algebra
of Section 7.

Let $\sF_n$ be the family of regular unital embeddings
\[
\phi : A(G_n) \otimes M_r \to A(G_n) \otimes M_s, ~~ \mbox 
{ for } r < s,
\]
which preserve the  2 by 2 block structure. Thus  $\phi$ may be
indicated  as
\[
\phi = \left[
\begin{array}{cc} \phi_1&\phi_{12}\\ 0&\phi_2 \end{array}\right]
\]
where $\phi_1, \phi_2 : \bbC^n \otimes M_r \to \bbC^n \otimes M_s$ are
C*-algebra maps and 
\[
\phi_{12}: M_{n,m} \otimes M_r \to M_{n,m} \otimes M_s
\]
is an appropriate bimodule map.
We shall show that $\sF_n$ is a not a functorial family.
The key construction for this is to obtain
an irregular factorisation of a regular embedding \\
$\theta : A(G_n) \to
A(G_n) \otimes M_{n^2}$ which is the direct sum
of $n^2$ maps $\theta_{ij}$ arising from automorphisms $\sigma_{ij}
= \sigma_1^i \times \sigma^j_2$ of $G_n$ where $\sigma_1 , \sigma_2$ are
 cyclic shifts  of the range and source vertices.
In Donsig and Power \cite{don-pow-2} this was obtained for $n = 2$ by a
seemingly fortuitous ad hoc argument.

Let $n\ge 2$ be an integer and let $w$ be a primitive root of unity. Let
$U$ be the $n \times n$ unitary matrix
\[
U = (u_{ij}) = (\frac{w^{(i-1)(j-1)}}{\sqrt{n}})
\]
and let $S$ be the cyclic forward shift in $M_n$ for the standard basis.

Let $(f_{ij})$ be the standard matrix unit system for $M_n$ where $M_n$
is viewed as the $(1,2)$ block subspace of $A(G_n)$.
Define 
$\phi : A(G_n)  \to A(G_n) \otimes M_n$ 
to be the restriction of the unique C*-algebra map $\tilde{\phi}$ between
the generated C*-algebras for which
\[
\tilde{\phi} : \left[\begin{array}{cc}
0& f_{ij}\\
0&0
\end{array}\right] \to  
\left[\begin{array}{cc} 0 &((S^*)^{i-1}US^{j-1}) \otimes 
e_{ij}\\ 0&0 \end{array}\right] 
\]
where $(e_{ij})$ is the standard matrix unit system for the tensor factor
$M_n$.
Since \[
\phi: \left[\begin{array}{cc}
0 & f_{ij}\\
0&0
\end{array}\right] \to 
(Ad Z) \circ 
\phi  \left( \left[\begin{array}{cc}
0 & U \otimes e_{ij}\\
0&0
\end{array}\right] \right)
\]
with $Z$ a $ 2 \times 2$ block diagonal unitary operator in $C^*(A(G_n)
\otimes M_n)$ it is clear that there is such a C*-algebra map. Note that
$\phi$ is not a regular embedding. For example, for $n = 3$ 
the  embeddings has multiplicity
3 and 
{\small
\[
\phi(f_{11}) = \frac{1}{\sqrt{3}}\left[\begin{array}{ccc|ccc|ccc}
1&0&0&1&0&0&1&0&0\\
0&0&0&0&0&0&0&0&0\\
0&0&0&0&0&0&0&0&0\\
\hline
1&0&0&w&0&0&w^2&0&0\\
0&0&0&0&0&0&0&0&0\\
0&0&0&0&0&0&0&0&0\\
\hline
1&0&0&w^2&0&0&w^4&0&0\\
0&0&0&0&0&0&0&0&0\\
0&0&0&0&0&0&0&0&0\\
\end{array}\right].
\]
}
Since the block entries have norm $\frac{1}{\sqrt{3}}$ the embedding
is not regular.

Let $\overline{U}$ be the complex congugate of the matrix $U$ and define
\[
\psi :  A(G_n) \otimes M_n \to (A(G_n) \otimes M_n)  \otimes M_n
\] to be the unique star algebra homomorphism such that for $x$ in the
tensor factor,
\[
\psi : \left[\begin{array}{cc} 0 &f_{ij}\\ 0&0 \end{array}\right]\otimes x
\to \left(\left[\begin{array}{cc} 0 & (S^*)^{i-1}US^{j-1}\\ 0&0 \end{array}\right] 
\otimes e_{ij}\right) \otimes x.
\]

\begin{lma}
The map $\psi \circ \phi$  is inner conjugate to
the 1-decomposable embedding $\theta$.
\end{lma}

\begin{proof}
Note that
\[
(S^*)^{i-1}US^{j-1} = (\frac{w^{(k+i-1)(l+j-1)}}{\sqrt{n}})^n_{k,l=1}
\]
and that
{\small
\[
(k+i-1)(l+j-1)-(s+k-1)(t+l-1) = (i-1)(j-1)-(s-1)(t-1)+k(j-t)+l(i-t).
\]}
\newpage

Thus we have the following  calculation.
\[
\begin{array}{ll}
(\psi \circ \phi)_{12}(f_{ij}) & = \psi_{12}(((S^*)^{i-1}US^{j-1})\otimes e_{ij})\\
& \\
& =  \psi_{12}({\sqrt{n}}^{-1}( \Sigma_{k,l} {w^{(k+i-1)(l+j-1)}}f_{k,l}) \otimes e_{ij})\\
& \\
& = {\sqrt{n}}^{-1} \Sigma_{k,l} {w^{(k+i-1)(l+j-1)}}(\psi_{12}(f_{k,l}\otimes e_{ij}))\\
&\\
& =  {\sqrt{n}}^{-1} \Sigma_{k,l} {w^{(k+i-1)(l+j-1)}}
{\sqrt{n}}^{-1}(\Sigma_{s,t} {w^{(s+k-1)(t+l-1)}}f_{st}\otimes e_{kl}\otimes e_{ij})\\
&\\
&= \Sigma_{s,t} f_{st} \otimes (\Sigma_{kl}w^{(i-1)(j-1)-(s-1)(t-1)}
(\frac{w^{j-t})^k(w^{i-s})^l}{n} e_{kl})\otimes e_{ij}\\
&\\
& = \Sigma_{s,t} f_{st} \otimes (X^{ij}_{st}) \otimes e_{ij}
\end{array}
\]
where $X^{ij}_{st}$ is a unimodular multiple of the rank one partial
isometry
\[
Y^{ij}_{st} = \Sigma_{k,l} \frac{(w^{j-t})^k(w^{i-s})^l}{n} e_{kl}.
\]
We now want to show that the composition $\psi \circ \phi$ is not merely
locally regular, which is what the above calculation shows, but that it
is regular, that is, 1-decomposable. (A
star extendible locally regular map
not be regular, as we have seen.)

Let $(g_1,g_2,\dots ,g_n)$ in $\bbC^n$ be the basis with
\[
g_i = (w^i, w^{2i}, \dots , w^{ni})/\sqrt{n},
\]
so that
\[
Y^{ij}_{st} = g_{j-t} \otimes \overline{g}_{i-s},
\] where 
$g_1 \otimes \overline{g}_2$ indicates the rank one operator for which 
$(g_1 \otimes \overline{g}_2)(h) = ~<h,g_2>g_1$.

Observe that the embedding $\eta : A(G_n) \to A(G_n) \otimes M_{n^2}$
for which
\[
\eta_{12}(f_{ij}) = \Sigma_{s,t} (f_{s,t} \otimes g_{j-t,i-s} \otimes
e_{ij})
\]
is a regular star extendible embedding unitarily equivalent to a map
$\theta$ determined by two cyclic shifts
as indicated above. Thus $\psi \circ \phi$ is the
composition $\eta \circ \theta$ where
\[
\eta : M_{2n} \to M_{2n} \otimes M_{n} \otimes M_{n}
\]
is the linear Schur product map given by
\[
\eta(f_{st}\otimes e_{kl}\otimes e_{ij} 
= w^{(i-1)(j-1)-(s-1)(t-1)} f_{st}\otimes e_{kl}\otimes e_{ij}.
\]

Although these unimodular coefficients do not form a cocycle, that is,
$\eta$ is not realisable as a diagonal unitary conjugation, the
restriction of $\eta$ to the span of the matrix units
\[
\{
f_{s,t} \otimes g_{j-t,i-s} \otimes
e_{ij}
 : 1 \le s,t \le n,   1 \le i,j \le n\}
\]
is a cocycle. This may be checked directly. Alternatively  note
that since $\psi \circ \phi = \eta \circ \theta$, the map 
$\psi \circ \phi$ is the orthogonal direct sum of $n^2$ star extendible embeddings
$\eta \circ \theta_{ij}$. Since $\psi \circ \phi$ is
star extendible so too is each map $\theta_{ij}$. That $\eta$ is diagonally
implementable on the ranges of the multiplicity one maps 
$\theta_{ij}$ follows from
the fact that an isometric Schur product map on the bipartite graph is
diagonally implementable.
(In fact this property is shown to hold for any digraph algebra in Davidson and
Power \cite{dav-pow}.)
Since there is a diagonal partition of the identity operator which
dominates the ranges of the maps $\eta \circ \theta_{ij}$ it follows that
$\eta \circ \theta$ is diagonally conjugate to $\theta$, as desired.
\end{proof}

Let $\sG_n \subseteq \sF_n$ denote the closed subfamily of regular maps
whose indecomposables are the multiplicity one embeddings corresponding
to the automorphisms of $G_n$.
The arguments above show that $\sF_n$ and $\sG_n$ are closed families
which do not satisfy the factorisation property. In fact these families
are not even functorial.

\begin{thm}
For $n = 2,3, \dots  $ the familes $ {\sF}_n$ and ${\sG}_n$ are not functorial.
\end{thm}

\begin{proof}
Let $\phi_k = \phi \otimes \id, \psi_k = \psi \otimes \id$, be the maps
\[
\phi \otimes \id : A(G_n) \otimes M_{n^2k} \to (A(G_n) \otimes M_n) \otimes
M_{n^2k},
\]
\[
\psi \otimes \id : A(G_n) \otimes M_{n^{2k+1}} \to (A(G_n) \otimes M_n) 
\otimes M_{n^{2k+1}}.
\]
Then for all $k$ the compositions $\psi_k \circ \phi_k$ are regular
by the last lemma. Also, since $\overline{w}$ is also a primitive root of
unity the lemma shows that the compositions 
$\phi_{k+1} \circ \phi_k$ are regular.
Thus the maps $\phi_k, \psi_k$ provide a commuting diagram between two
regular systems in $\Sys (\sG_n)$ consisting of irregular maps.
Moreover it is clear that the crossover maps of any subdiagram are
necessarily irregular, and so the theorem is established.
\end{proof}

\noindent {\bf Subalgebra positions.}

In \cite{pow-matroid}, 
 we analysed irregular factorisations in the case of the
4-cycle algebra $A(G_2)$ and showed  that there is a converse to the
construction above in the following sense.
If $\Phi$ is an irregular star-extendible isomorphism between the
systems ${\cal A}, {\cal A}'$ in $\Sys ({\cal G}_2)$ then necessarily 
 ${\cal A}$ and $ {\cal A}'$ are systems determined by compositions of
embeddings of type $\theta$. In particularly ${\cal A}$ and $ {\cal A}'$
are regularly isomorphic by some other isomorphism $\Psi$. Because of
this ${\cal G}_2$ is an isoclassic family. Noting that ${\cal G}_2$ is
the family of rigid embeddings we deduce that the equivalence between
(i),(ii) and (iv) in Theorem 7.4 also holds for the case $n=2$ of
4-cycle algebras.

Let ${\cal G}_2^{UHF} \subseteq {\cal G}_2$ be the subfamily of unital
systems ${\cal A}$ for which the algebraic direct limit has the form
\[
A_0 = \left[
\begin{array}{cc}D_0&M_0\\
&D_0
\end{array}\right]
\]
where $D_0$ is a unital ultramatricial algebra.
It was shown in \cite{pow-matroid}
 how the inclusion
\[
D_0 \oplus D_0 \subseteq B_0 = C^*(A_0)
\]
 determines the pair $A_0,
A_0^*$ and  therefore how  the $K_0H_1$ classification scheme of Donsig
and Power \cite{don-pow-2} gives a classification 
scheme for  these positions.
In themselves each summand $D_0$ has Jones index 2 in the 
corresponding corner algebra of $B_0$ and these positions are unique,
in analogy with (although more elementary than)
 Goldman's theorem for index 2 subfactors.
 Thus the invariants may be viewed 
as determining
 the relative position of index 2 subcorners in the superalgebra.
By taking weak closures in the tracial representation one obtains
unital inclusions $R \oplus R \subseteq R$, where $R$ is the
hyperfinite II$_1$ factor  (with common Cartan masa)
and where,  again, the summands have index 2 in the 
corners. In this case,
as one would expect, almost all of the $K_0H_1$ invariants evaporate.
The residue turns out to be 
 the $H_0H_1$ coupling invariant for partial isometry
homology. One can reinterpret  the weak
closures in the language of subfactor theory 
in terms of an $R-R$-bimodule ${}_RM^{\alpha}_R$ determined by
a symmetry $\alpha$, in which case
 the $H_0H_1$ coupling invariant corresponds to Connes spectral
invariant \cite{con} for the symmetry.

We expect that similar techniques will show that for the other
bipartite graphs the family ${\cal G}_n$ is isoclassic and 
hence that one can similarly obtain complete
dimension module  invariants for the bipartite (algebraic)
limit algebras with respect to these regular embeddings. Once
again this leads to invariants for regular subalgebra positions 
$D_0 \oplus D_0 \subseteq M_0$ (of higher Jones index)
and connections with subfactor theory.
However, there are more possibilities for irregular factorisations of
1-indecomposable embeddings and it becomes interesting
 to determine the appropriate subfactor setting.
This connection should lead to information on the number of
approximately inner
equivalence classes of standard diagonals
(cf Donsig Power \cite{don-pow-1}) in the bipartite limit
algebras (for ${\cal G}_n$) and may  perhaps shed light on the
longstanding problem of the automorphic uniqueness of standard
diagonals in regular limit algebras.

\end{spacing}
 \end{document}